\documentclass[preprint,12pt]{elsarticle}

\usepackage{graphicx}
\usepackage{epsfig}
\usepackage{epstopdf}
\usepackage{psfrag}
\usepackage{pgfplots}
\usepackage{subcaption}

\usepackage{amssymb}

\usepackage{amsthm}
\usepackage{amsmath}

\usepackage{lineno}
\usepackage{bbold}
\usepackage{lipsum}
\usepackage{amsfonts}
\usepackage{dsfont}
\usepackage{epstopdf}
\usepackage{algorithmic}
\usepackage{cleveref}
\usepackage{xcolor}
\crefformat{equation}{(#2#1#3)}

\usepackage{url}

\usepackage{tikz,tkz-euclide}
\usetikzlibrary{shapes.misc}
\usetikzlibrary{calc}
\tikzset{cross/.style={cross out, draw=black, minimum size=2*(#1-\pgflinewidth), inner sep=0pt, outer sep=0pt},
cross/.default={1mm}}

\newtheorem{proposition}{Proposition}

\newtheorem{lemma}{Lemma}

\newtheorem{remark}{Remark}

\journal{Journal of Computational Physics}

\begin{document}

\newcommand{\Surf}{{\partial \Omega}}
\newcommand{\Vol}{{\Omega}}
\newcommand{\dVol}{{\partial\Omega}}
\newcommand\txtred[1]{{\color{red}#1}}
\newcommand\txtblue[1]{{\color{blue}#1}}
\newcommand{\x}{x}
\newcommand{\y}{y}
\newcommand{\half}{\frac{1}{2}}
\newcommand{\BigO}[1]{\ensuremath{\operatorname{O}\bigl(#1\bigr)}}
\newcommand{\Order}{\mathcal{O}}
\newcommand{\U}{{\mathbf{U}}}
\newcommand{\Ue}{{\mathbf{V}}}
\newcommand{\Ux}{{\mathbf{U}_x}}
\newcommand{\Vv}{{\mathbf{V}}}
\newcommand{\vv}{{\overline{v}}}
\newcommand{\G}{{\mathbf{G}}}
\newcommand{\A}{\mathbf{A}}
\newcommand{\Pp}{\mathcal{P}}
\newcommand{\Ppo}{{\mathcal{\overline{P}}}}
\newcommand{\Q}{\mathcal{Q}}
\newcommand{\Qx}{\mathcal{Q}_{x}}
\newcommand{\Qxe}{\mathcal{Q}_x^e}
\newcommand{\Qxx}{\mathcal{Q}_{xx}}
\newcommand{\Qxi}{\mathcal{Q}_{x^i}}
\newcommand{\B}{\mathcal{B}}
\newcommand{\Lagr}{\mathcal{L}}
\newcommand{\lagr}{\mathcal{l}}
\newcommand{\Dx}{\mathcal{D}_{x}}
\newcommand{\Dxi}{\mathcal{D}_{x^i}}
\newcommand{\Dxa}{\mathcal{D}_{x^1}}
\newcommand{\Dxb}{\mathcal{D}_{x^2}}
\newcommand{\Dxc}{\mathcal{D}_{x^3}}
\newcommand{\Bx}{\B\Dx}
\newcommand{\Lx}{\mathbb{L}_x}
\newcommand{\Di}{\mathcal{D}_{AD}}
\newcommand{\Dii}{\mathcal{D}_{AD_i}}
\newcommand{\De}{\mathcal{D}^{e}}
\newcommand{\Dxx}{\mathcal{D}_{xx}}
\newcommand{\Ee}{\mathcal{E}}
\newcommand{\R}{\mathcal{R}}
\newcommand{\nn}{\nonumber}
\newcommand{\Nn}{\ell}
\newcommand{\Nns}{\ell^\star}
\newcommand{\Mm}{\mathcal{M}}
\newcommand\norm[1]{\left\lVert#1\right\rVert}
\newcommand\bint[3]{\left. #1 \right|_{#2}^{#3}}
\newcommand\sg[1]{\sigma^#1}
\newcommand\eps{\varepsilon}
\newcommand\epst{\overline{\varepsilon}}
\newcommand{\V}{\overline{u}}
\newcommand{\F}{\mathcal{F}}
\newcommand{\Aa}{\mathcal{A}}
\newcommand\bb[1]{\textcolor{blue}{#1}}
\newcommand\rred[1]{\textcolor{red}{#1}}
\newcommand\ggreen[1]{\textcolor{green}{#1}}

\begin{frontmatter}



\title{An SBP-SAT Continuous Galerkin Finite Element Formulation for Smooth and Discontinuous Fields}


\author{A.G. Malan$^\dagger$}
\author{J. Nordstr$\ddot{\textnormal{o}}$m\corref{cor1}$^{\ddagger\star}$}
\ead{jan.nordstrom@liu.se$^\star$}
\cortext[cor1]{Corresponding Author, Dept. of Mathematics, Linköping Univ., Sweden.}
\address{$^\dagger$InCFD Research Group,  Dept. of Mechanical Eng., Univ. of Cape Town, South Africa \\ $^\ddagger$Dept. of Mathematics, Applied Mathematics, Linköping Univ., Sweden, \\ $^\star$Dept. of Mathematics and Applied Mathematics, Univ. of Johannesburg, South Africa} 

\begin{abstract}
The high-order accurate continuous Galerkin finite element method offers attractive computational efficiency for computational fluid dynamics. A challenge is however spurious oscillations which result for convection dominated flows over discontinuities.  To derive a continuous Galerkin scheme for both smooth and discontinuous fields we start by first writing the scheme in Summation-by-Parts (SBP) form for a single element mesh.  Boundary conditions are applied weakly via Simultaneous-Approximation-Terms (SAT) and Gauss-Labotto quadrature employed in the interest of computational efficiency.  We then show that the stable single element baseline scheme in SBP-SAT form extends trivially to a provably stable multi-element formulation.  Next, we develop provably stable element based Galerkin-weighted artificial dissipation operators to deal with spurious oscillations over shocks while retaining high order accuracy elsewhere.  The resulting scheme achieves super-convergence with accuracy of $\Order(p+2)$ when using $p^{th}$ order Lagrange polynomials for smooth fields.  The developed dissipation operators furnish WENO like behaviour over discontinuities while retaining high order accuracy elsewhere for both linear and non-linear wave propagation problems.

\end{abstract}


\end{frontmatter}

\section{Introduction}
\label{SpatDisc}

After decades of advances in computer processing speed, computational cost remains a bottleneck in computational fluid dynamics (CFD). This is especially acute in cases involving phenomenon which travel long distances or oscillate over extended time periods.  A ground breaking procedure demonstrating the ability of high-order accurate methods to model this was presented in \cite{Kreiss1972}.  Since then high-order accurate methods have matured for certain industrial strength applications \cite{Peraire2017,Jameson2018,Changfoot19,Ilangakoon2020,Sevilla2021, Ilangakoon2022,Nordstrom2023}.  These examples are representative of the main competing high-order accurate methods in CFD today i.e. the discontinuous Galerkin (DG) finite element, flux reconstruction (FR) and finite difference (FD) methods.  The high-order accurate continuous Galerkin (CG)  finite element method (FEM) has to date received relatively little attention despite being known to be highly efficient \cite{Huerte2013}.  A key reason for this is its lack of inherent numerical dissipation and resulting spurious oscillations when modelling convection dominated flows,  particularly for  fields containing discontinuities.

The \textit{Summation-By-Parts} (SBP) framework was originally developed for the FEM \cite{Kreiss1974,Kreiss1977}.  Since then, it has been mainly developed for high-order accurate FD methods.  Provable stability is achieved with ease when combining SBP with weakly imposed \textit{Simultaneous-Approximation-Term} (SAT) boundary conditions proposed in \cite{Carpenter1994} and further developed in \cite{Carpenter1991,Nordstrom1999,Nordstrom2001}.  The SBP-SAT methodology (reviewed in \cite{Svard2014,Fernandez2020}) has recently been  applied to the CG method  for both linear \cite{2020JanCG1} and non-linear \cite{2023JanCG2} advection of smooth fields.  We build on this work by developing an SBP-SAT compliant provably stable artificial dissipation (AD) scheme for the Galerkin weighted FEM.  This to deal with the transport of discontinuous fields with-in elements when using high-order Lagrange polynomials ($>2^{nd}$ order) for CG and DG schemes.

The difficulty in removing spurious oscillations over discontinuities is a property that CG and DG share with vertex centered finite volume (FV),  central difference FD and spectral methods.  In the case of FV this has historically been addressed by upwinding \cite{Changfoot19} or added artificial dissipation (AD) \cite{Jameson1995,Nordstrom2006_2,Pattinson2007}. AD is applied to reduce the order of accuracy over the shock while retaining full accuracy elsewhere.  In high-order accurate FD using the SBP-SAT framework, spurious oscillatory behaviour has been removed via provably stable AD and filtering procedures \cite{Mattsson2003,Lundquist2020}. The fundamental development in \cite{Mattsson2003} was later mimicked in  \cite{Diener2007,Svard2009,Zingg2018,Ranocha2018}, and we apply it here in the Galerkin weighted setting. Spectral element methods have also been augmented using AD \cite{Tadmor1989,Maday1993} but WENO like behaviour over discontinuities not achieved. Finally,  the DG method was recently extended to include CG based AD to damp oscillations interior to elements \cite{Peraire2023},  but stability was not proven rigorously. 

In contrast to other methods,  the use of AD in the CG method has received comparatively little attention \cite{Hughes1986,Manzari1998} and is limited to low order Lagrange polynomials without provable stability. We aim to address this.  First,  by building on \cite{2020JanCG1,2023JanCG2} and following \cite{Nordstrom2017} we derive a boundary condition imposition procedure for the scalar advection-diffusion equation and prove it energy bounded in its \textit{continuous form}.  Second, we mimic the continuous analysis in semi-discrete Galerkin weighted form for a single element using SBP-SAT principles to arrive at an energy stable base-line scheme.  In the third step,  we show that the single element formulation in SBP-SAT form automatically yields a stable multi-element high-order CG formulation. Our focus on coupling of elements via CG is due to a wealth of literature already available for the coupling in DG \cite{Hesthaven1996,Gassner2013,Kopriva2021}.  Following \cite{Mattsson2003,Ranocha2018}, we in the fourth and final step construct an element based AD formulation in SBP form and prove it stable. The AD essentially reduces the spatial order in an \textit{element wise} manner applied only in the vicinity of the discontinuity while retaining high-order accuracy elsewhere.  The resulting AD operator is naturally also applicable to the DG method. 

For computational efficiency we employ Gauss-Lobatto (GL) quadrature,  which achieves super convergence in accuracy of $\Order(p+2$) ($p$ being the Lagrange polynomial order) for the linear smooth advection-diffusion equation.  In addition, we demonstrate that WENO like behaviour is achieved over discontinuities for both linear and non-linear hyperbolic  problems when using Lagrange polynomials up to fourth order.  In summary we develop a provably stable high-order accurate CG scheme on SBP-SAT form for smooth and discontinuous fields. As a side effect we obtain an element based AD operator which is directly applicable to DG.

The remainder of the paper is organised as follows.  The next section prepares for the SBP framework by deriving energy estimates using integration-by-parts.  This is followed by a semi-discrete SBP-SAT treatment on a single element, and subsequently multi-element stability is proven for the CG method.  The section to follow details the element based AD scheme suitable for CG and DG.  The performance of the resulting CG method with the new AD is finally assessed via several numerical examples.

\section{The General Framework in Continuous Form}

We consider the 1D scalar advection-diffusion equation in the domain $x \in \Vol=[x_0,x_1]=[0,1]$ as follows:
\begin{eqnarray}
u_t + a u_x & = \left( \eps u_{x} \right)_x,&~~t\geq 0 
\label{eq:1D-adv-diff} \nn \\
G_0 (u_0(t),\eps_0(t)) & =  g_0(t), &~~x =0  \\
G_1 (u_1(t),\eps_1(t)) & =  g_1(t), &~~x=1  \nn \\ 
u(x,0) & =  f(x),&~~x \in \Vol. \nn
\end{eqnarray}

\noindent In (\ref{eq:1D-adv-diff}), $a, \eps \geq 0$ and $a$ is constant while $\eps=\eps(x,t)$. Further, $G_0$ and $G_1$ are boundary condition (BC) operators at $x_0$ and $x_1$ respectively while the boundary and initial data are  respectively $g_{0},g_{1}$ and $f$.  We employ a Robin BC at the left inflow boundary and a Neumann condition at the right outflow boundary as

\begin{equation}\label{eq:RobinBCs}
G_0(u,\eps) = \bint{(a u - \eps u_{x})}{0}{} = g_0;~~G_1(u,\eps) =\bint{- \eps u_{x}}{1}{}=g_1,
\end{equation}

The BCs will be implemented weakly via \textit{simultaneous approximation terms} (SAT) \cite{Carpenter1994} in the numerical setting.  Following \cite{Nordstrom2017},  we simulate this in the continuous setting such that the governing equation  reads

\begin{equation}
u_t + a u_x = \left( \eps u_{x} \right)_x  + \mathbb{SAT}.
\label{eq:1D-adv-diff-SAT}
\end{equation}

\noindent Here we,  with a slight abuse of notation,  introduce a SAT like penalty term

\begin{equation} \label{eq:1D_cont SAT1}
\mathbb{SAT}  =  L_0  \left( \sigma_0 \left[ (a u - \eps u_{x})-g_0 \right] \right) +  L_1 \left(\sigma_1 \left[ (- \eps u_{x})-g_1 \right] \right)
\end{equation}

\noindent where $\sigma_0$ and $\sigma_1$ are penalty parameters to be determined. The boundary conditions are added weakly using the \textit{lifting operators} $L_0,L_1$ \cite{Brezzi01} defined as

\begin{equation*}
\int_\Vol \psi L_i(\phi) d\Vol = \int_{\partial \Vol_i} \psi \phi dS.
\end{equation*}

\noindent Here $\psi$ and $\phi$ are smooth scalar functions and $\partial \Vol_i$ denotes the bounding surface of $\Vol$ at $x_i$. 

\begin{remark}
The variable coefficient $\eps$ will inspire our development of AD.
\end{remark}

\subsection{The Continuous Energy Estimate}

The energy method is applied by multiplying the governing equation (\ref{eq:1D-adv-diff-SAT}) with the solution, integrating by parts over $\Vol$ followed by rearranging to yield

\begin{align}
\frac{1}{2}\frac{d}{dt} \norm{u}_\Vol^2  +  \norm{ \sqrt{\eps} u_x}_\Vol^2 = BT
\label{eq:adv_diffenergy0.2}
\end{align}
\noindent where $\norm{\phi}_\Vol^2 = \int_\Vol \phi^2 dx$ for some scalar function $\phi$. Further 
\begin{align}
BT& = \bint{u\left( \frac{a}{2}u-\eps u_x \right)}{1}{0} +  \int_{\Vol} u \mathbb{SAT} ~d\Vol. \nn \\
    & = - \frac{a}{2} \left(u_0^2 +u_1^2 \right) + {\color{red} u_0 \left(a u_0 -\eps_0 u_{x_0} \right)} +{\color{red} u_1 \left(\eps_1 u_{x_1} \right)} +  \int_{\Vol} u \mathbb{SAT} ~d\Vol. \nn
\end{align}
	
\noindent The terms in red above may be positive and must be cancelled via the SAT boundary conditions.  The lifting operators give
\begin{equation}\label{eq:adv_diffenergy0.3}
 \int_{\Vol} u \mathbb{SAT} ~d\Vol = \sigma_0 u_0 \left[ \bint{\left( a u-\eps u_x \right)}{0}{}- g_0 \right] + u_1 \sigma_1 \left[- \bint{\eps u_x}{1}{} - g_1 \right].
\end{equation}

\noindent An energy estimate exists if the resulting $BT \leq 0$ for $g_0=g_1=0$ \cite{Nordstrom05}.  This will be the case provided
\begin{equation}\label{eq:cont_penalties_0}
\sigma_0 = -1 ;~~\sigma_1 = 1.
\end{equation}

\noindent The parameter values (\ref{eq:cont_penalties_0}) in combination with time integration yields the energy estimate
\begin{equation} \label{eq:cont_energy}
\norm{u}_\Vol^2 + \int_0^t 2 \norm{ \sqrt{\eps} u_x}_\Vol^2  dt \leq \norm{f}^2_\Vol
\end{equation}

By mimicking the above continuous procedure using the SBP-SAT framework \cite{Carpenter1991}, we now seek to prove that the 1D Galerkin FEM with a similar boundary condition (BC) imposition is stable.  We aim for a semi-discrete  estimate of type (\ref{eq:cont_energy}).

\begin{figure}[ht]
\centering
\def\svgwidth{0.7\textwidth}
\begingroup%
  \makeatletter%
  \providecommand\color[2][]{%
    \errmessage{(Inkscape) Color is used for the text in Inkscape, but the package 'color.sty' is not loaded}%
    \renewcommand\color[2][]{}%
  }%
  \providecommand\transparent[1]{%
    \errmessage{(Inkscape) Transparency is used (non-zero) for the text in Inkscape, but the package 'transparent.sty' is not loaded}%
    \renewcommand\transparent[1]{}%
  }%
  \providecommand\rotatebox[2]{#2}%
  \newcommand*\fsize{\dimexpr\f@size pt\relax}%
  \newcommand*\lineheight[1]{\fontsize{\fsize}{#1\fsize}\selectfont}%
  \ifx\svgwidth\undefined%
    \setlength{\unitlength}{370.4669141bp}%
    \ifx\svgscale\undefined%
      \relax%
    \else%
      \setlength{\unitlength}{\unitlength * \real{\svgscale}}%
    \fi%
  \else%
    \setlength{\unitlength}{\svgwidth}%
  \fi%
  \global\let\svgwidth\undefined%
  \global\let\svgscale\undefined%
  \makeatother%
  \begin{picture}(1,0.72925752)%
    \lineheight{1}%
    \setlength\tabcolsep{0pt}%
    \put(0,0){\includegraphics[width=\unitlength,page=1]{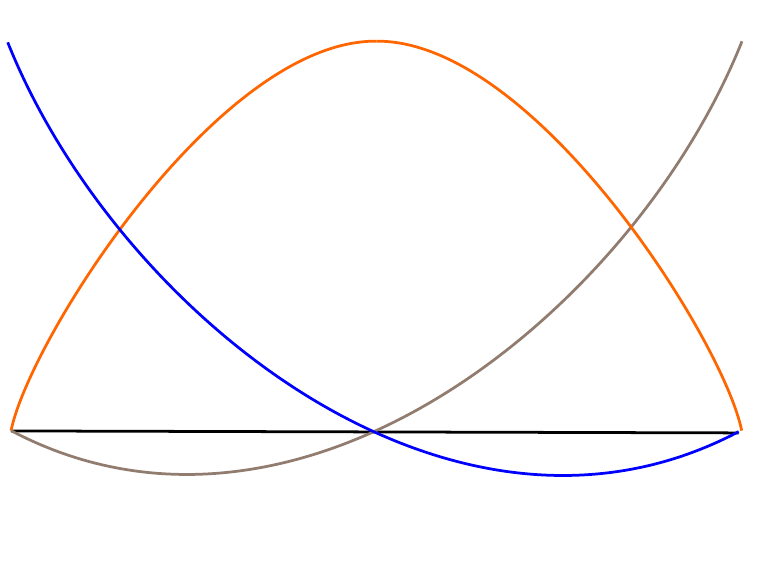}}%
    \put(0.00509961,0.0977173){\makebox(0,0)[lt]{\lineheight{1.25}\smash{\begin{tabular}[t]{l}$0$\end{tabular}}}}%
    \put(0.06047679,0.62321067){\makebox(0,0)[lt]{\lineheight{1.25}\smash{\begin{tabular}[t]{l}$ \Nn_0$\end{tabular}}}}%
    \put(0.47447988,0.70669834){\makebox(0,0)[lt]{\lineheight{1.25}\smash{\begin{tabular}[t]{l}$\Nn_1$\end{tabular}}}}%
    \put(0.83970561,0.60121798){\makebox(0,0)[lt]{\lineheight{1.25}\smash{\begin{tabular}[t]{l}$\Nn_2$\end{tabular}}}}%
    \put(0.46846358,0.0977173){\makebox(0,0)[lt]{\lineheight{1.25}\smash{\begin{tabular}[t]{l}$1$\end{tabular}}}}%
    \put(0.94576227,0.0977173){\makebox(0,0)[lt]{\lineheight{1.25}\smash{\begin{tabular}[t]{l}$2$\end{tabular}}}}%
    \put(0,0){\includegraphics[width=\unitlength,page=2]{1D-elements_CG.pdf}}%
  \end{picture}%
\endgroup%

	\caption{1D quadratic ($p=2$) element where $\Nn_i$ are the Lagrange basis polynomials.  }
	\label{fig:FEM_1D}
\end{figure}

\section{Single Element CG Framework on SBP-SAT Form}

We will initially  consider a 1D domain $\Vol$ consisting of a single element. The mesh contains $N+1$ nodes numbered from $0$ to $N$. The trial solution over the element is given by
\begin{equation}
\V(x,t) = \Lagr^T (x) {\U(t)}
\label{eq:lagr1D}
\end{equation}
\noindent where $\Lagr(x)^T = \left[\Nn_{0}(x), \Nn_{1}(x),..., \Nn_{p}(x)\right]$ consists of $p+1$ Lagrange polynomials $\Nn_{i}(x)$ of order $p$. Here
\begin{equation*}
\Nn_i(x)= \frac{(x-x_{0})...(x-x_{i-1})(x-x_{i+1})...(x-x_{p})}{(x_i-x_{0})...(x_i-x_{i-1})(x_i-x_{i+1})...(x_i-x_{p})}.
\end{equation*}
\noindent where $x_i$ are $p+1$ Gauss-Lobatto (GL) points which coincide with the mesh nodes.  Therefore $p=N$.  Further,  $\U(t)=~[\V_0, \V_1, ... , \V_p]^T$ consists of $p+1$ time-dependant coefficients. 

\noindent Note that $\V = u$ if the exact solution $u$ is a polynomial of order $p$ or less.  The approximations $u_x (x,t) \approx \Lagr_{x}^{T}(x) \U(t)$ and $u_{xx}(x,t) \approx \Lagr_{xx}^{T}(x) \U(t)$ follow directly from (\ref{eq:lagr1D}).  Here the derivative of the vector of Lagrange polynomials evaluated at node $i$, which is located at $x_i$, is  $\Lagr^T_{x_i}= \left[ \Nn_{x_0}(x_i), \Nn_{x_1}(x_i), . . . , \Nn_{x_p}(x_i) \right]$ where  $\Nn_{x_j}(x) = d \Nn_j/dx$.

We next integrate the governing equation (\ref{eq:1D-adv-diff-SAT}) over $\Vol$ together with the Galerkin weighting and introduce  ansatz (\ref{eq:lagr1D}). This gives
\begin{align}
\frac{\partial}{\partial t}  \int_{\Vol} \Lagr \Lagr^T dx \U & + a  \int_{\Vol} \Lagr \Lagr_{x}^{T} dx \U 
= \int_{\Vol} \Lagr \left( \eps \Lagr_{x}^{T} \right)_x dx \U +  \int_{\Vol} \Lagr \mathds{SAT} dx,
\label{eqn:adv_diff_discr_1D}
\end{align}

\noindent where $\int_{\Vol}$ denotes the spatial integral over the element.  Further,
\begin{align} \label{eq:SAT_1D_4}
\int_{\Vol} \Lagr {\mathds{SAT}} dx & =   \sigma_0 \Lagr_0  \left[ \left( a\Lagr^T_0 - \eps_0 \Lagr^T_{x_0} \right)\U - \Lagr^T _0  \mathbf{g} \right] \nn \\ 
& + \sigma_N \Lagr_N \left[- \eps_N \Lagr^T_{x_N} \U - \Lagr^T_{N}  \mathbf{g} \right] 
\end{align}
\noindent where  $\mathbf{g} = [g_0, 0, ... ,0, g_N]^T$ and $\sigma_{0,N}$ are respectively the left and right penalty coefficients which correspond to (\ref{eq:cont_penalties_0}).

We will employ GL quadrature for computing the spatial integrals as
\begin{equation}
\int_{\Vol} f(x) dx = \sum_{i = 0}^{p} w_i f(x_i).
\label{eq:gauss-q}
\end{equation}

\noindent Here $x_i$ and $w_i$ respectively denote the GL node coordinate and corresponding quadrature weight. (\ref{eq:gauss-q}) is exact for polynomials $f(x)$ of order $\leq (2p - 1)$.

\subsection{The Discrete Operators on SBP-SAT Form}
We now develop the spatial operators in (\ref{eqn:adv_diff_discr_1D}) such that it reads
\begin{equation}
\Pp \U_t + a  \Qx  \U =  \Qxx \U + \int_{\Vol} \Lagr {\mathds{SAT}} d\Vol. 
\label{eq:adv_diff_disc_2}
\end{equation}

\noindent In (\ref{eq:adv_diff_disc_2}),  $\Pp$ is the diagonal mass matrix and $\Qx$ is the weak form first derivative operator. Further, $\Qxx=\Qxx\left(\eps\right)$ is the weak form second derivative operator, in this case a function of $\eps(x,t)$.

We start the development of (\ref{eq:adv_diff_disc_2}) with the mass matrix,  and get
\begin{align}\label{eq:mass_1D}
\Pp \approx & \int_{\Vol} \Lagr \Lagr^T dx  \approx \sum_{i =0}^p w_i  \Lagr(x_i) \Lagr^T(x_i)=  diag[w_0,w_1,...,w_p]. 
\end{align}

\noindent $\Qx$ is next cast into SBP form by applying integration-by-parts to the advective term \cite{Nordstrom2006} followed by GL quadrature:
\begin{align}
\Qx = \int_{\Vol} \Lagr \Lagr_x^T dx = \bint{\Lagr \Lagr^T}{0}{N} - \int_{\Vol} \Lagr_x \Lagr^T dx =  \mathcal{B} - \Qx^{T}.
\label{eq:Qe_1D}
\end{align}

\noindent Here $\Qx$ is almost skew-symmetric since $\mathcal{B} = \bint{\Lagr \Lagr^T}{0}{N} = (\Qx + \Qx^T)=diag[ -1, 0,... ,0,1]$. Importantly, $\B$ is an operator, non-zero only on the domain boundary.

We similarly cast the diffusion term into SBP form by applying integration-by-parts:
\begin{eqnarray} \label{eq:discr_diff_1D_0}
\Qxx(\eps) = \int_{\Vol}\mathcal{L} \left( \eps    \Lagr_{x}^{T} \right)_x dx & =  &  \bint{\Lagr \eps \Lagr_x^T}{0}{N} - \int_{\Vol} \mathcal{L}_x   \eps \Lagr_{x}^{T} dx.
\end{eqnarray}

\noindent To develop the discrete operators for (\ref{eq:discr_diff_1D_0}), we digress to define the required expressions.  For this we will need the following notations:
\begin{equation*}
\mathcal{E} = diag[\eps_0, \eps_1, ...,\eps_p];~~\sqrt{\mathcal{E}} = diag[\sqrt{\eps_0}, \sqrt{\eps_1}, ...,\sqrt{\eps_p}].
\end{equation*}

\begin{proposition} By employing GL quadrature, the weak form $2^{nd}$ second derivative operator (\ref{eq:discr_diff_1D_0}) becomes
\begin{equation}\label{eq:prop2_3}
\Qxx(\eps) = \mathcal{E} \Bx -  \left(\sqrt{\mathcal{E}} \Dx\right)^T\Pp \left(\sqrt{\mathcal{E}} \Dx\right).
\end{equation}
\noindent where $\Dx = \Pp^{-1} \Qx$ is the strong form first derivative operator.
\begin{proof} 

Consider the general scalar $v(x)$ and vector $\Vv$ such that $v(x) = \Lagr^T \Vv$.  Applying GL quadrature to the weak form $1^{st}$ derivative operator yields
\begin{equation*}\label{eq:prop2_eq0}
\Qx \Vv = \int_{\Vol}  \Lagr \Lagr^T_x  dx \Vv = \Pp \mathbb{L}_x^T \Vv  \textnormal{~~where~~}
\mathbb{L}_x^T= \begin{bmatrix} \Nn_{x_0}(x_0) &  \cdots & \Nn_{x_p}(x_0)  \\ 
\vdots & \ddots & \vdots \\
 \Nn_{x_0}(x_p)  & \cdots &  \Nn_{x_p}(x_p) \\
\end{bmatrix}. 
\end{equation*}

\noindent The discrete first derivative is therefore computed as $\Dx \Vv =  \Pp^{-1} \Qx \Vv =\Lx^T \Vv$. This clarifies the left identity of the right-hand-side in (\ref{eq:discr_diff_1D_0}) and (\ref{eq:prop2_3}). Using this we consider the last term on the right-hand-side of (\ref{eq:discr_diff_1D_0}):
\begin{align*}
~~~~~~~~~~~~~~~~~\Vv^T \int_{\Vol} \Lagr_x \eps \Lagr_x^T dx  \Vv & = \int_{\Vol} \left( \sqrt{\eps} \Lagr_x^T\Vv \right)^T  \left( \sqrt{\eps} \Lagr_x^T\Vv \right)  dx \\
& = \left(\sqrt{\mathcal{E}} \Lx^T \Vv \right)^T\Pp \left(\sqrt{\mathcal{E}} \Lx^T\Vv \right) \\
& =	 \Vv^T\left(\sqrt{\mathcal{E}} \Dx \right)^T\Pp \left(\sqrt{\mathcal{E}} \Dx\right) \Vv. ~~~~~~~~~~~~~~~~~~~~~~~~~~\square
\end{align*}
\phantom\qedhere \end{proof}\end{proposition}

\subsection{The Discrete Energy Estimate}
All differential operators equired for (\ref{eq:adv_diff_disc_2}) have now been developed on SBP-SAT form. Next we  mimic the continuous energy analysis in discrete form by multiplying (\ref{eq:adv_diff_disc_2}) with $\U^T$ followed by applying Proposition 1:
\begin{equation}
 \U^T \Pp \U_t +  \left( \sqrt{\mathcal{E}} \Ux \right)^T \Pp \left( \sqrt{\mathcal{E}} \Ux \right) + a \U^T \Qx  \U =  \U^T \mathcal{E} \Bx \U + \U^T\int_{\Vol} \Lagr {\mathds{SAT}} dx.  
\label{eq:adv_diff_disc_energy1}
\end{equation}

\noindent The last term in (\ref{eq:adv_diff_disc_energy1}) is
\begin{align*}
\U^T \int_{\Vol} \Lagr {\mathds{SAT}} d\Vol  
= &  \U^T  \left\{ \sigma_0 \Lagr_0  \left[ \left( a\Lagr^T_0 - \eps_0 \Lagr^T_{x_0} \right)\U - \Lagr^T _0  \mathbf{g} \right] \right.  \nn \\ 
& ~~~~\left. + \sigma_N \Lagr_N \left[- \eps_N \Lagr^T_{x_N} \U - \Lagr^T_{N}  \mathbf{g} \right] \right\} \nn \\
= & \sigma_0 \V_0 \left[ a \V_0-\eps_0 \V_{x}(x_0) - g_0 \right] + \V_N \sigma_N \left[-\eps_N \V_{x}(x_N) - g_N \right].
\label{eq:1D_boun_fin}
\end{align*}

\noindent We next seek $\sigma_0$ and $\sigma_N$ to ensure a bounded energy norm and thus stability.

\begin{proposition}
$\sigma_{0}=-1$ and $\sigma_{N}=1$ in (\ref{eq:adv_diff_disc_energy1}) leads to an energy estimate.

\begin{proof}
By taking the transpose of (\ref{eq:adv_diff_disc_energy1}) and  adding it to itself yields
\begin{equation}
 \frac{d}{dt} \norm{\U}^2_{\Pp} + 2 \norm{\sqrt{\mathcal{E}} \U_{x}}^2_{\Pp} = \mathbb{BT}.
\label{eq:adv_diff_norm}
\end{equation}

\noindent Here $\norm{\U}^2_{\Pp} = \U^T \Pp \U$, $\norm{\sqrt{\mathcal{E}} \U_{x}}^2_{\Pp} = \left( \sqrt{\mathcal{E}} \Ux \right)^T \Pp \left( \sqrt{\mathcal{E}} \Ux \right)$ and
\begin{align}
\mathbb{BT} = & -a \U^T \mathcal{B} \U  + 2 \U^T \mathcal{E} \Bx \U + 2~\U^T \int_{\Vol} \Lagr {\mathds{SAT}} dx \\
          = &-a\left( \V_0^2+\V_N^2 \right) + 2 (\sigma_0+1) \left[ a \V_0 \V_0-\eps_0 \V_0 \V_{x}(x_0) \right] + \nn \\
          & 2(\sigma_N-1) \left[-\eps_N \V_N \V_{x}(x_N)  \right]. \nn
\end{align}

\noindent where we have used zero boundary data. For stability, $\mathbb{BT}$ must be non-positive which requires
\begin{equation*}
\sigma_0 = -1; ~~\sigma_N =1.
\end{equation*}

\noindent Inserting this result and integrating (\ref{eq:adv_diff_norm}) over time:
\begin{equation}
\norm{\U}^2_{\Pp}  + 2 \int_0^t \norm{\sqrt{\mathcal{E}} \U_{x}}^2_{\Pp} dt \leq \norm{\mathcal{\overline{F}}}_\Pp^2
\label{eq:adv_diff_norm2}
\end{equation}
\noindent where $\mathcal{\overline{F}}$ is the discrete initial data.  Note that (\ref{eq:adv_diff_norm2}) mimics (\ref{eq:cont_energy}). \end{proof} \end{proposition}

\begin{remark} Consistent with the continuous form (\ref{eq:adv_diffenergy0.2}), the single element Galerkin weighted scheme on SBP-SAT form has no energy growth terms internal to the domain. The only growth terms are due to the boundaries. The SAT boundary condition removes these and is therefore in effect a damping term.  No added dissipation is required for stability \cite{2020JanCG1} even for $\eps=0$.
\end{remark}
\begin{remark} The proof of Proposition 2 mimics the continuous energy analysis,  including the choice of penalty parameters.  The SBP-SAT procedure guarantees a numerical treatment consistent with the continuous one. Once the continuous procedure for an energy bound is known, stability follows almost automatically by applying a similar procedure to the discrete form \cite{Nordstrom2017}.
\end{remark}

\section{Multi-Element CG Framework on SBP-SAT Form}\label{sec:mult-elem}

Up to this point we developed all operators on a single element mesh.  In the CG multi-element setting,  element based operators are merged.  We exemplify the procedure by considering a case consisting of a left (L) and right (R) element and nodes numbered sequentially from $0$ to $N=(2p+1)$.  A global operator $\Aa$ is obtained by merging the two single element operators as follows: 
\begin{align*}
\Aa & = \sum_e \Aa^e = \txtred{\Aa^L} + \txtblue{\Aa^R} \nn \\ 
& = \begin{bmatrix}
\txtred{\Aa_{00}^L} & \txtred{\cdots} & \txtred{\Aa_{0p}^L} &  & 0 \\
\txtred{\vdots} & \txtred{\ddots} & \txtred{\vdots} &  &  \\
\txtred{\Aa_{p0}^L} & \txtred{\ddots} & (\txtred{\Aa_{pp}^L}+\txtblue{\Aa_{00}^R}) & \txtblue{\cdots} & \txtblue{\Aa_{0p}^R} \\
&  & \txtblue{\vdots} & \txtblue{\ddots} & \txtblue{\vdots} \\
0 &  & \txtblue{\Aa_{p0}^R} & \txtblue{\cdots} & \txtblue{\Aa_{pp}^R}
\end{bmatrix}.
\end{align*}
\noindent Here the subscripts denote the row and column of the respective single element operators.  Using this merging procedure, all internal (non-boundary) diagonals of the global $\Qx$ are zero due to the symmetry properties of the single element operator $\Qxe$ where $\Q_{x_{pp}}^e=-\Q_{x_{00}}^e$. This is illustrated by considering the two element mesh and merging two quadratic elements:
\begin{align*}
\Qx & = \sum_e \Qxe = \txtred{\Qx^L} + \txtblue{\Qx^R} \nn \\ 
& = \frac{1}{6} \begin{bmatrix} \begin{array}{ccccc} 
\txtred{-3} & \txtred{4} & \txtred{-1} & 0  & 0 \\
\txtred{-4} & \txtred{0}  & \txtred{4} & 0 & 0 \\
\txtred{1} & \txtred{-4} & (\txtred{3}-\txtblue{3})  & \txtblue{4} & \txtred{-1} \\
0 & 0 & \txtblue{-4} & \txtblue{0}  & \txtblue{4} \\
 0 & 0 & \txtblue{1} & \txtblue{-4} & \txtblue{3}
\end{array} \end{bmatrix}.
\end{align*}
 
\noindent The globally merged $\Qx$ therefore remains almost skew symmetric when the transpose is added (see (\ref{eq:Qe_1D})) i.e. $\Qx+\Qx^T = \B$.  The symmetry properties are therefore retained automatically when extending the single element formulation to the CG multi-element setting.

The inter-element contributions of the non-symmetric $\Qxx(\eps)$ terms similarly do not appear due the $\Ee \B \Dx$ term  in  (\ref	{eq:prop2_3}) only being non-zero on domain boundaries. This achieves flux continuity at element boundaries which is consistent with the continuous equations (this is achieved in standard CG \cite{Zien2000} by explicitly removing $\bint{\Lagr \eps \Lagr_x^T}{0}{p}$ in (\ref{eq:discr_diff_1D_0}) at inter-element boundaries).

Therefore, the CG scheme written in SBP-SAT form retains the symmetry properties of the single element case.  We summarise this as:
\begin{proposition} The symmetry properties of the single element operators $\Qx$ and $\Qxx$ are retained in the multi-element CG setting if written in SBP form.
\end{proposition}

We finally prove that the last term on the right-hand-side of (\ref{eq:prop2_3}) remains dissipative in the CG multi-element setting. To this end we need Lemma 1.

\begin{lemma}
Consider a two element mesh consisting of a left (L) and a right (R) element.  Let the single element operators $\A^{L,R}$ be positive semi-definite i.e. $\Vv^T \Aa^{L,R} \Vv \geq 0$ where $\Vv$ is an arbitrary vector.  Then the merged operator $\Aa=\Aa^L + \Aa^R$ is also positive semi-definite.

\begin{proof} 
In the merged two-element setting we get
\begin{align*}
\Vv^T \sum_e \Aa^e \Vv & = \begin{bmatrix} \vv_0 \\ \vdots \\ \vv_p \\ \vdots \\ \vv_{N} \end{bmatrix}^T
\begin{bmatrix}
\Aa_{00}^L & \cdots & \Aa_{0p}^L &  & 0 \\
\vdots & \ddots & \vdots &  &  \\
\Aa_{p0}^L & \cdots & (\Aa_{pp}^L+\Aa_{00}^R) & \cdots & \Aa_{0p}^R \\
&  & \vdots & \ddots & \vdots \\
0 &  & \Aa_{p0}^R & \cdots & \Aa_{pp}^R
\end{bmatrix}
\begin{bmatrix} \V_0 \\ \vdots \\ \vv_p \\ \vdots \\ \vv_{N} \end{bmatrix} \\
&= \begin{bmatrix} \vv_0 \\ \vdots \\ \vv_{p} \end{bmatrix}^T
 \Aa^L
 \begin{bmatrix} \vv_0 \\ \vdots \\ \vv_{p} \end{bmatrix} + 
 \begin{bmatrix} \vv_p \\ \vdots \\ \vv_{N} \end{bmatrix}^T
 \Aa^R
 \begin{bmatrix} \vv_p \\ \vdots \\ \vv_{N} \end{bmatrix} \geq 0.~~~~~~~~~~~~~~~~~\square
\end{align*} \phantom\qedhere \end{proof}\end{lemma}

For $\Aa^e=\left(\sqrt{\Ee} \Dx^e \right)^T \Pp^e \left( \sqrt{\Ee} \Dx^e \right)$, Lemma 1 states that the last term on the right-hand-side of (\ref{eq:prop2_3}) remains positive semi-definite in the multi-element setting. This in turn implies that the energy norm in the multi-element setting satisfies

\begin{equation}
\sum_e \norm{\U}^2_{\Pp^e}  + 2 \int_0^t \sum_e
\norm{\sqrt{\Ee} \U_{x}}^2_{\Pp^e} dt \leq \sum_e \norm{\mathcal{\overline{F}}}_{\Pp^e}^2.
\label{eq:adv_diff_norm3}
\end{equation}

We have now proven the following proposition.

\begin{proposition}
The stability of the single-element CG scheme in SBP-SAT form (\ref{eq:adv_diff_disc_energy1})  is retained in the multi-element setting subsequent to merging.
\end{proposition}

\section{Galerkin Weighted Artificial Dissipation}

The stable CG scheme derived above is next extended by adding artificial dissipation (AD).  Our development of the Galerkin weighted AD operator is inspired by the variable diffusion term treatment in the base-line scheme.

\subsection{Single Element AD}
The single element based advection equation with AD, ignoring BCs is
\begin{equation} \label{eq:adv}
\Pp \U_t + a {\Qx} \U + a \Di \U = 0.
\end{equation}

\noindent The dissipation operator $\Di$ with $p$ components is defined as
\begin{equation}\label{eq:AD_f}
\Di = \sum_{i=1}^p \gamma_i \Dii = \sum_{i}^{p} \gamma_i \int_{\Vol^e} \Lagr_{x^i} \eps_i \Lagr_{x^i}^T dx
\end{equation}

\noindent where $i$ denotes the order of the Lagrange polynomial derivative e.g.  $\Lagr_{x^3}=\Lagr_{xxx}$.  Further,  $\eps_i \geq 0$ and $\gamma_i$ can attain the values $[-1,0,1]$. The single component dissipation term in (\ref{eq:AD_f}) is computed as
\begin{equation}\label{eq:diss_00}
\Dii = \left(\sqrt{\Ee_i} \Dxi^T \right)^T \Pp \left(\sqrt{\Ee_i} \Dxi^T \right)
\end{equation}
\noindent (\ref{eq:diss_00}) is positive semi-definite for $\Ee_i \geq 0$, just as the diffusion term in the base-line scheme. This definiteness property of the single element AD is trivially retained in the multi-component case (\ref{eq:AD_f}) for $\gamma_i \geq 0$ due Lemma 1. 

In addition we can in the special case of a constant $\eps_i$ over an element prove that we can allow for some $\gamma_i=-1$. In this case we consider

\begin{equation}
\Di = \sum_{i=1}^p \epst_i \Dxi \Pp \Dxi^T
\end{equation}

\noindent where $\epst_i = \gamma_i \eps_i$ and constant over an element.  Consider  as an example the quadratic element  $(p=2)$ for which

\begin{align} \label{eq:Di2}
\Di & = \epst_1 \Dxa \Pp \Dxa^T + \epst_2 \Dxb \Pp \Dxb^T \nonumber \\
&  = \frac{1}{6} \begin{bmatrix}
7 \epst_1 + 12 \epst_2 & -8 \epst_1 -24 \epst_2 & \epst_1 + 12 \epst_2 \\
-8 \epst_1 - 24 \epst_2 & 16 \epst_1 + 48 \epst_2 & -8 \epst_1 - 24 \epst_2 \\
\epst_1 + 12 \epst_2 & -8 \epst_1 -24 \epst_2 & 7 \epst_1 + 12 \epst_2 
\end{bmatrix}.
\end{align}

\noindent Note that while $\epst_1$ is always positive,  $\epst_2$ is allowed to be negative.

\begin{proposition} The AD operator (\ref{eq:Di2}) is positive semi-definite provided that
\begin{equation} \label{eq:Diss2_Stability}
\epst_1 > 0~~~and~~~\epst_2 \geq -\frac{1}{3} \epst_1
\end{equation}
\noindent where $\epst_1$ and $\epst_2$ are constants.
\begin{proof}  Positive semi-definiteness is guaranteed provided that the eigenvalues of (\ref{eq:Di2}) are all non-negative and real.  By rotating the $\Di^e$ matrix to diagonal form via a non-singular rotation matrix $R$ we find $R^T \Di^e R=\Lambda $ where
\begin{equation}
\Lambda  = diag \left[ 0,  \epst_1,  4 \epst_1 + 12 \epst_2 \right].
\end{equation}

\noindent By Sylvester's law \cite{Horn1990} we conclude that all eigenvalues are non-negative provided that (\ref{eq:Diss2_Stability}) holds. \end{proof}\end{proposition}
  
We next consider a cubic element ($p=3$) with spatially constant $\epst_i$ in which case the added dissipation term would be
\begin{equation} \label{eq:Di3}
\Di =  \epst_1 \Dxa \Pp \Dxa^T + \epst_2 \Dxb \Pp \Dxb^T + \epst_3 \Dxc \Pp \Dxc^T
\end{equation}
\begin{proposition} The AD operator (\ref{eq:Di3}) is positive semi-definite provided that
\begin{equation}  \label{eq:Diss3_Stability}
\epst_1 > 0;~~~\epst_{2} \geq -\frac{\epst_1}{3};~~~\epst_3 \geq  -\frac{\epst_1}{45}-\frac{\epst_{2}}{3}
\end{equation}
\noindent  where $\epst_1$, $\epst_2$ and $\epst_3$ are constants.
\begin{proof}  An energy norm is guaranteed provided that the eigenvalues of (\ref{eq:Di3}) are all non-negative and real.  The  entries of the rotated diagonal matrix are
\begin{equation} \label{eq:eigen3}
\Lambda  = diag \begin{bmatrix}  
0 \\ 50 \epst_1/3 + 50 \epst_2 \\
\F + \sqrt{\F - 100\epst_1 (45\epst_3 + \epst_1 + 15\epst_2)} \\ 
\F - \sqrt{\F - 100\epst_1 (45\epst_3 + \epst_1 + 15\epst_2)} \\ 
\end{bmatrix}
\end{equation}

\noindent where $\F = 17\epst_1 + 225\epst_2 + 675\epst_3$.

\noindent The second eigenvalue in (\ref{eq:eigen3}) is positive for $\epst_1 \geq 0$ if $\epst_2 \geq -\epst_1/3$. Considering the remaining eigenvalues and noting that $\epst_1 \geq 0$, we require 
\begin{align}\label{eq:eig3_2}
45\epst_3 + \epst_1 + 15\epst_2 \geq 0;~\F \geq 0~~and~~\F - 100\epst_1 (45\epst_3 + \epst_1 + 15\epst_2) \geq 0
\end{align}

\noindent The first condition implies
\begin{align} \label{eq:eig3_3}
\epst_3 = -\frac{\epst_1}{45} - \frac{\epst_2}{3} + \delta~~where~~\delta \geq 0
\end{align}

\noindent By inserting (\ref{eq:eig3_3}) into the second requirement in (\ref{eq:eig3_2}),  we find it automatically satisfied since $\F \geq 0$ if $2\epst_1 + 675\delta \geq 0$. Next, (\ref{eq:eig3_3}) is inserted into the final requirement in (\ref{eq:eig3_2}) to yield the requirement:
\begin{align} \label{eq:eig3_4}
f(\delta) = 4\epst_1^2 - 1800\epst_1\delta + \delta^2455625 \geq 0
\end{align}

\noindent We first note that $f(\delta=0) = 4 \epst_1^2 \geq 0$. Next we seek the minimum $f(\delta)$ since $d^2 f/d\delta^2 = 911250 > 0$. $f(\delta)$ will be minimum where $df/d \delta = 911250 \delta - 1800 \epst_1 = 0$. Therefore $f_{min}(\delta=4\epst_1/2025) = 20/9~\epst_1^2$ and hence the final condition in (\ref{eq:eig3_2}) holds. \end{proof} \end{proposition}

\begin{remark} The procedure for definiteness in Propositions 5 and 6 can be extended to higher order AD, but requires even more complex algebra.
\end{remark}
\begin{remark}The choice of constant $\epst_i$ over the element is only required if we allow for negative coefficients.  We have not explored this option in the numerical examples to follow but include it as a stable  alternative.  If all $\epst_i \geq 0$ we can of course allow for variable coefficients with-in the element.
\end{remark}
\begin{remark}
A variable dissipation coefficient can also be used to damp oscillations locally at discontinuities.  This can be done by letting $\epst_i$ be non-zero only in the vicinity of the discontinuity. 	In the numerical examples to follow we use $\epst_i$ as a step function,  being non-zero only around the discontinuity.
\end{remark}

\subsection{Multi-Element AD}
For multiple elements we have the following proposition.

\begin{proposition}
In the multi-element $CG$ setting where each element based $\Di^e$ is positive definite, the merged dissipation operator is stable.
\begin{proof}
See Lemma 1. \end{proof}
\end{proposition}

\begin{remark}
For DG the multi-element AD coupling is trivially stable since elements are not merged.
\end{remark}

\subsection{A Few Instructive AD Examples}
The effect of the added dissipation on a multi-element mesh is exemplified by first considering linear elements ($p=1$) and assuming $\epst_1 \geq 0$ as well as $a>0$. The element based first derivative and dissipation operators are
\begin{equation} \label{eq:DisP1}
\Qx =  \half
\begin{bmatrix} -1 & 1  \\ -1 & 1 \end{bmatrix};~
\Di = \int_{\Vol^{e}} \Lagr_{x} \epst_1  \Lagr_{x}^T dx = \frac{\epst_1}{2} \begin{bmatrix} 1 & -1  \\ -1 & 1 \end{bmatrix};
\end{equation}

\begin{figure}[ht]
\centering
\def\svgwidth{1.1\textwidth}
	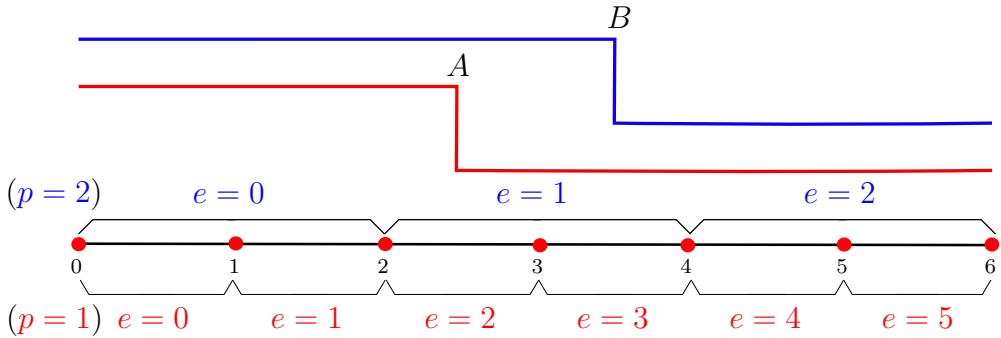
	\caption{Schematic of a $7$ node mesh consisting of $6$ linear (\txtred{$p=1$}) or $3$ quadratic (\txtblue{$p=2$}) elements.  Shocks A and B are respectively depicted in red and blue. }
	\label{fig:FEM_1D_L1}
\end{figure}

\noindent Consider the $6$ element \txtred{$p=1$} mesh with the discontinuous field $A$ shown in Figure~\ref{fig:FEM_1D_L1}.  In this case:
\begin{equation}\label{eq:Qx1}
\Qx = \half \begin{bmatrix} \begin{array}{ccc|cccc} 
 -1 & 1  & 0 & 0 & 0  & 0 & 0 \\
 -1 & 0  & 1 & 0 & 0  & 0  & 0\\
 0  & -1 & 0 & 1 & 0  & 0  & 0\\
0 & 0  & -1 & 0 & 1  &  0  & 0\\
0 & 0 & 0  & -1 & 0 &  1& 0 \\
0 & 0 & 0  & 0 & -1 &  0& 1 \\
0 & 0 & 0  & 0  & 0 & -1 & 1 \\
 \end{array} \end{bmatrix}
\end{equation}

\noindent where the vertical line marks the shock position $A$.  We let $\epst_1>0$ over element $\txtred{e=2}$ and add it to the assembled $\Qx$ to yield
\begin{equation}
 \Qx + \Di = \half \begin{bmatrix} \begin{array}{cccc|cccc} 
\frac{row}{column} & \text{\scriptsize{0}} & \text{\scriptsize{1}} & \text{\scriptsize{2}} & \text{\scriptsize{3}} & \text{\scriptsize{4}} & \text{\scriptsize{5}} & \text{\scriptsize{6}} \\
\text{\scriptsize{0}} &  -1 & 1 & 0 & 0  & 0 & 0 & 0 \\
\text{\scriptsize{1}} & -1 & 0  & 1 & 0 & 0 & 0 & 0 \\
\text{\scriptsize{2}} &0  & -1 & \txtred{\epst_1}  & 1-\txtred{\epst_1} & 0 & 0 & 0 \\
\text{\scriptsize{3}} &0 & 0 & -1 \txtred{-\epst_1} & \txtred{\epst_1}& 1& 0 & 0  \\
\text{\scriptsize{4}} &0 & 0 & 0  & -1 &0 & 1 & 0 \\
\text{\scriptsize{5}} &0 & 0 & 0  & 0 & -1  & 0 & -1 \\
\text{\scriptsize{6}} &0 & 0 & 0  & 0 & 0  & -1 & -1 \\
 \end{array} \end{bmatrix}.
\end{equation}

\noindent The dissipation terms are coloured red for  clarity.  Next set $\epst_1 = 1$ to yield
\begin{equation}\label{eq:Q1_sc}
 \Qx + \Di = \half \begin{bmatrix} \begin{array}{cccc|cccc} 
\frac{row}{column} & \text{\scriptsize{0}} & \text{\scriptsize{1}} & \text{\scriptsize{2}} & \text{\scriptsize{3}} & \text{\scriptsize{4}} & \text{\scriptsize{5}} & \text{\scriptsize{6}} \\
\text{\scriptsize{0}} & -1 & 1 & 0 & 0  & 0 & 0 & 0 \\
\text{\scriptsize{1}} & -1 & 0  & 1 & 0 & 0 & 0 & 0 \\
\text{\scriptsize{2}} &0 & -1 & \txtred{1}  & \txtred{0} & 0 & 0 & 0 \\
\text{\scriptsize{3}} &0 & 0  & \txtred{-2} & \txtred{1}  & 1  & 0 & 0 \\
\text{\scriptsize{4}} &0 & 0 & 0  & -1 & 0 & 1& 0 \\
\text{\scriptsize{5}} &0 & 0 & 0  & 0 & -1  & 0 & -1 \\
\text{\scriptsize{6}} & 0 & 0 & 0  & 0 & 0  & -1 & -1 \\
 \end{array} \end{bmatrix}.
\end{equation}

\noindent which results in an upwind biassed stencil in rows $2$ and $3$.   This is consistent with what is done in FV and second order accurate FD.  Note that since $\epst_1 \geq 0$ the eigenvalues of $\Di$ are non-negative and stability is preserved.

We next consider  a quadratic element, for which
\begin{align} \label{eq:Q2}
\Qxe & = \frac{1}{6} \begin{bmatrix}  -3 & 4 & -1 \\ -4 & 0 & 4 \\ 1 & -4 & 3 \end{bmatrix}.
\end{align}

\noindent For a constant $\epst_i$ on each element, the dissipation $\Di^e$ is given by (\ref{eq:Di2}). Consider now the discontinuous field $B$ in Figure~\ref{fig:FEM_1D_L1} in conjunction with the quadratic (\txtblue{$p=2$}) element mesh.  By way of example,  with $\epst_1=\epst_2=0$ on $\txtblue{e=0}$, and $\epst_1 = 1/3$; $\epst_2 = 1/18$ on elements \txtblue{$e=1$} and \txtblue{$e=2$}, we obtain the following upwind-biased and provably stable operator
 \begin{equation} \label{eq:quad_v2}
 \Qx +  \Di = \frac{1}{6} \begin{bmatrix} \begin{array}{ccccc|ccc} 
\frac{row}{column} & \text{\scriptsize{0}} & \text{\scriptsize{1}} & \text{\scriptsize{2}} & \text{\scriptsize{3}} & \text{\scriptsize{4}} & \text{\scriptsize{5}} & \text{\scriptsize{6}} \\
 \text{\scriptsize{0}} &-3 & 4 & -1 & 0  & 0 & 0 & 0 \\
\text{\scriptsize{1}} & -4 & 0  & 4 & 0 & 0 & 0 & 0 \\
\text{\scriptsize{2}} &1 & -4 & \txtred{3}  & \txtred{0} & \txtred{0} & 0 & 0 \\
\text{\scriptsize{3}} &0 & 0  & \txtred{-8} & \txtred{8}  &\txtred{0} & 0 & 0 \\
\text{\scriptsize{4}} &0 & 0 & \txtred{2} & \txtred{-8}  & \txtred{6} & \txtred{0} & \txtred{0}\\
\text{\scriptsize{5}} &0 & 0 & 0  & 0 & \txtred{-8}  & \txtred{8} & \txtred{0} \\
\text{\scriptsize{6}} & 0 & 0 & 0  & 0 &  \txtred{2}  &  \txtred{-8} &  \txtred{6} \\
 \end{array} \end{bmatrix}.
\end{equation}

\noindent The vertical line shows the shock position. The entries in the coefficient matrix (\ref{eq:quad_v2}) which are affected by the dissipation are in red (rows $3$ to $6$) and are now upwind biased.

\subsection{Variable Length Elements}\label{sec:paraCG}

Elements of varying lengths (standard length is required for efficient GL quadrature) will be employed for the numerical examples to follow. In this case, the discretization of the advection-diffusion equation is transformed to
\begin{align}
\frac{\partial}{\partial t} \int_{\Vol^e} \Lagr \Lagr^T d\xi \U J & + a \int_{\Vol^e} \Lagr \Lagr_{\xi}^{T} d\xi \U + a  \Di \U \nn \\ & = \int_{\Vol^e} \Lagr \left( \eps \Lagr_{\xi}^{T} \right)_\xi  d\xi \U J^{-1} +  \int_{\Vol^e} \Lagr \mathds{SAT} d \xi J.
\label{eqn:adv_diff_discr_1D_para}
\end{align}

\noindent In (\ref{eqn:adv_diff_discr_1D_para}) the global ($x$) and transformed coordinates ($\xi$) are related as follows:
\begin{equation*}
x(\xi) = \frac{1-\xi}{2}x_0 + \frac{\xi+1}{2}x_p~~~\text{and}~~~J = \frac{dx}{d\xi}=\frac{x_p-x_0}{2}.
\end{equation*}

\noindent The standard coordinate system spans $-1 \leq \xi \leq 1$ with $x_{0,p}$ as the left most and right most global nodal coordinates of the element.

Next we write the AD in terms of the values derived to ensure stability:
\begin{equation}\label{eq:diss_para}
\Di = \sum_{i=1}^p  \int_{\Vol^{e}} \Lagr_{\xi^i} \left( \frac{\epst_i^\prime}{J^{[2i-1]}} \right)  \Lagr_{\xi^i}^T d\xi = \sum_{i=1}^p  \int_{\Vol^{e}} \Lagr_{\xi^i} \epst_i \Lagr_{\xi^i}^T d\xi.
\end{equation}

\noindent where $ \epst_i^\prime$ is the dissipation coefficient in the transformed coordinate frame.  The choice of $\epst_i$ is invariant with respect to element size.  This is possible since the key purpose of AD is to reduce the accuracy of the weak form first derivative which remains invariant with respect to element size i.e. $\Qx(x) = \Qx(\xi)$. 

The required AD in the transformed coordinate system automatically reduces  with reduction in element size. This is evident by relating the dissipation in the transformed coordinate frame $\epst_1^\prime,\epst_2^\prime,\epst_3^\prime$ to $\epst_1,\epst_2,\epst_3$ as follows:
\begin{equation*}
\epst_1^\prime = \epst_1 J;~~~\epst_2^\prime = \epst_2 J^3;~~~\epst_3^\prime = \epst_3 J^5.
\end{equation*}
\noindent Clearly less dissipation is obtained as the mesh resolution increases since the Jacobian $(J)$ reduces in magnitude.

\begin{figure}[ht!]
\centering
\begin{tikzpicture}[font=\footnotesize]
\begin{axis}[
    xmin = 0.0, xmax = 1,
    ymin = -0.1, ymax = 1.2,
    xtick distance = 0.2,
    ytick distance = 0.2,
    legend pos = south west,
    legend cell align={left},
    legend style={font=\footnotesize},    
    grid = both,
    minor tick num = 2,
    major grid style = {lightgray},
    minor grid style = {lightgray!25},
    width = 0.8\textwidth,
    height = 0.6 \textwidth,
    xlabel = {$x$},
    ylabel = {$u$},]
 

\addplot[
    thick,
    black,
    solid
] file[] {AdvDiff_10_P3_V.csv};
 
\addplot[
    only marks,
mark=otimes,
    blue,
] file[] {AdvDiff_10_P1_U_nodiss_10node.csv};

\addplot[
mark=square*,
    only marks,
    red,
] file[] {AdvDiff_10_P3_U_nodiss_10node.csv};

\legend{
   Exact, P1~(10~nodes),  P3~(10~nodes)}
\end{axis}
\end{tikzpicture}

\begin{tikzpicture}[font=\footnotesize]
\begin{axis}[
    xmin = 0.0, xmax = 1,
    ymin = -0.1, ymax = 1.2,
    xtick distance = 0.2,
    ytick distance = 0.2,
    legend pos = south west,
    legend cell align={left},
    grid = both,
    minor tick num = 2,
    major grid style = {lightgray},
    minor grid style = {lightgray!25},
    width = 0.8\textwidth,
    height = 0.6\textwidth,
    xlabel = {$x$},
    ylabel = {$u$},]
 

\addplot[
    thick,
    black,
    solid
] file[] {AdvDiff_40_P3_V.csv}; 

\addplot[
mark=square*,
    blue,
    solid
] file[] {AdvDiff_40_P4_U_nodiss_9node.csv};

\addplot[
mark=square,
    blue,
    solid
] file[] {AdvDiff_40_P4_U_diss_9node.csv};

\addplot[
    only marks,
mark=diamond*,
    red,
] file[] {AdvDiff_40_P4_U_nodiss_41node.csv};

\legend{
   Exact, p4~(9~node),  P4~(9~nodes) $\epst_1=\frac{3}{1000}$ $\epst_2=\frac{3}{250,000}$, P4~(41~nodes)}
\end{axis}
\end{tikzpicture}

\caption{Smooth advection-diffusion:  (Top) Exact solution for $a/\eps=10$ compared to solutions for $p=2$ and $p=3$.  (Bottom) Exact solution for $a/\eps=40$ compared to $p=4$.  Artificial dissipation is added only  to the $9$ node P4 mesh in red.}
\label{fig:advdiff_1}
\end{figure}
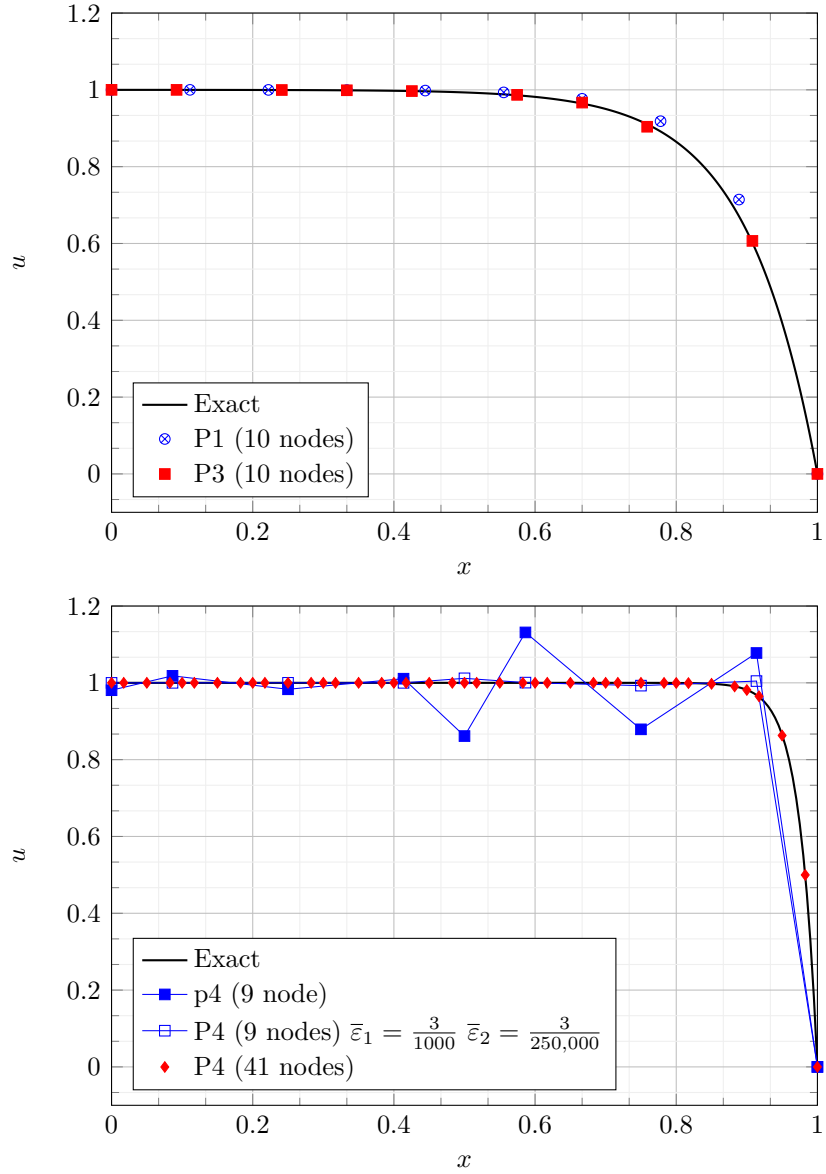

\begin{remark}
The reduction of Galerkin weighted AD as the Jacobian decreases proves that the developed AD is indeed artificial. This decrease with element size corresponds to those in FD \cite{Mattsson2003} and FV \cite{Nordstrom2006_2}.
\end{remark}

\section{Numerical Examples}

The  CG scheme including the AD is next applied to numerical examples varying from smooth to non-smooth and linear to non-linear. For this we employ Lagrange polynomials of order $p=1$ up to $p=4$ (denoted $P1$ to $P4$ in the figures). 

\begin{table}
\begin{small}
\begin{center}
\begin{tabular}{||r || c c | c c | c c | c c ||}
\hline \hline
& \multicolumn{8}{|c||}{\textbf{Advection-Diffusion: Smooth Field}}\\
\hline
\bf{Nodes} & \multicolumn{2}{|c|}{\bf{P1}} & \multicolumn{2}{|c|}{\bf{P2}} & \multicolumn{2}{|c|}{\bf{P3}} & \multicolumn{2}{|c||}{\bf{P4}} \\
     & $\epsilon_p$ & $\Order(\epsilon_p)$ & $\epsilon_p$ & $\Order(\epsilon_p)$  & $\epsilon_p(\epsilon_p)$ & $\Order(\epsilon_p)$ & $\epsilon_p$ & $\Order(\epsilon_p)$\\
\hline
& \multicolumn{8}{|c||}{$a/\eps=10$}\\
\hline
10  &  1.76E-02  &     &  4.44E-03  &     &  3.05E-03  &     &  2.47E-03  &   \\
19  &  4.17E-03  &  2.08  &  5.30E-04  &  3.61  &  1.56E-04  &  4.29  &  8.55E-05  &  4.85 \\
40  &  8.71E-04  &  2.02  &  2.39E-05  &  3.88  &  3.89E-06  &  4.78  &  4.92E-07  &  5.63 \\
85  &  1.87E-04  &  2.01  &  1.25E-06  &  3.97  &  8.74E-08  &  4.95  &  6.09E-09  &  5.92 \\
181  &  4.07E-05  &  2.00  &  5.97E-08  &  3.99  &  1.95E-09  &  4.99  &  6.38E-11  &  5.98 \\
361  &  1.01E-05  &  2.00  &  3.73E-09  &  4.00  &  6.11E-11  &  5.00  &  9.99E-13  &  6.00 \\
\hline
& \multicolumn{8}{|c||}{$a/\eps=40$}\\
\hline
10  &  1.4E-01  &     &  9.0E-02  &     &  8.6E-02  &     &  8.7E-02   &  \\
19  &  3.8E-02  &  1.88  &  2.4E-02  &  2.21  &  1.6E-02  &  2.39  &  1.6E-02  &  2.45 \\
40  &  7.4E-03  &  2.11  &  2.2E-03  &  3.00  &  1.1E-03  &  3.48  &  4.5E-04  &  3.89 \\
85  &  1.5E-03  &  2.07  &  1.5E-04  &  3.65  &  3.8E-05  &  4.39  &  9.8E-06  &  5.16 \\
181  &  3.3E-04  &  2.02  &  7.5E-06  &  3.91  &  9.6E-07  &  4.83  &  1.2E-07  &  5.74 \\
361  &  8.1E-05  &  2.00  &  4.8E-07  &  3.98  &  3.1E-08  &  4.96  &  2.0E-09  &  5.93 \\
\hline \hline
\end{tabular}
\caption{\label{tab:ad1} Smooth steady advection-diffusion CG related mesh convergence results for  $p=1$ to $p=4$.  } 
\end{center}
\end{small}
\end{table} 

\subsection{Linear Smooth Propagation}

In the first numerical test we consider the accuracy of the scheme for the steady advection-diffusion equation
\begin{align}
 a u_x &= \eps u_{xx};~~~\Vol \in [0,1];~~~u(0)  =  1;~~~u(1)  =  0\nn 
\end{align}

\noindent with the analytical solution
\begin{align*}
 u(x) &= 1 - \frac{e^{xa/\eps}-1}{e^{ a/\eps}-1} 
\end{align*}

\noindent We employ the in- and outflow BCs (\ref{eq:RobinBCs}) for the numerical solution and compute the resulting discretization errors via the P-Norm:
\begin{equation*}
\epsilon_p = \sqrt{(\U^T - \Ue^T) \Pp (\U - \Ue) }.
\end{equation*}

\noindent Here $\Ue$ is the exact solution  (the analytical solution inserted at nodes). The SAT penalties $\sigma_0 = -1 ;~\sigma_N = 1$ are employed throughout. The order of mesh convergence between successive meshes is computed as
 \begin{equation*}
\Order(\epsilon_p) = \frac{\log_{10}{\epsilon_{p_{2}}}-\log_{10}{\epsilon_{p_{1}}}}{\log_{10}{N_{1}}-\log_{10}{N_{2}}}
\end{equation*}

\noindent where $N$ is the number of nodes and the subscripts $1$ and $2$ denote two meshes with $N_2 > N_1$.

The mesh convergence study covers the cases of $a/\eps=10$ and $a/\eps=40$,  for which the exact solutions are compared to the computed values in Figure~\ref{fig:advdiff_1}.  As shown (top) there is a visible improvement in accuracy for the $a/\eps=10$ case if using $p=3$ as opposed to $p=1$ on a $10$ node mesh. Also depicted are the spurious oscillations resulting from solving the $a/\eps=40$ case using $p=4$ on a $9$ node mesh. These oscillations are due to insufficient mesh resolution (as discussed in \cite{Frenander2018}) and not due to an instability.  These may be reduced by either refining the mesh or via the addition of AD \cite{Malan2002} as shown in the bottom figure.

The order of accuracy achieved for a range of meshes is recorded in Table~\ref{tab:ad1}.  As shown,  super convergence of order $(p+2)$ is achieved throughout for $p>1$. This is the first article (to the best of our knowledge) to demonstrate such accuracy for CG via the application of weak boundary conditions. 

\begin{figure}[ht!]
\centering
\begin{tikzpicture}[font=\footnotesize]

\begin{axis}[
    xmin = 0.0, xmax = 1,
    ymin = 0.3333, ymax = 2.2,
    xtick distance = 0.25,
    ytick distance = 0.5,
    legend pos = north east,
    legend cell align={left},
    legend style={font=\footnotesize},
    grid = both,
    minor tick num = 2,
    major grid style = {lightgray},
    minor grid style = {lightgray!25},
    width = 0.8\textwidth,
    height = 0.6\textwidth,
    xlabel = {$x$},
    ylabel = {$u$},]
 

\addplot[
    thick, dashed,
    black,
] file[] {Adv_P3_V_init.csv}; 

\addplot[
thick,
    black,
] file[] {Adv_P3_V_end.csv}; 

\addplot[
only marks,
mark=triangle*,
    blue
] file[] {Adv_P1_U.csv}; 

\addplot[
only marks,
mark=square,
    blue
] file[] {Adv_P2_U.csv}; 

\addplot[
only marks,
mark=o,
    red
] file[] {Adv_P3_U_80_alpha_beta_gamma.csv}; 

\addplot[
mark=diamond*,
    black,
    dashed
] file[] {Adv_P4_U_80_alpha.csv}; 

\legend{Initial,Exact,P1, P2, P3, P4}

\end{axis}
\end{tikzpicture}

\begin{tikzpicture}[font=\footnotesize]
 
\begin{axis}[
    xmin = 0.0, xmax = 1,
    ymin = -0.22, ymax = 0.37,
    xtick distance = 0.1,
    ytick distance = 0.1,
    legend pos = north west,
    legend cell align={left},
    grid = both,
    minor tick num = 2,
    major grid style = {lightgray},
    minor grid style = {lightgray!25},
    width = 0.8 \textwidth,
    height = 0.6 \textwidth,
    xlabel = {$x$},
    ylabel = {$Error$},]
 

\addplot[
only marks,
mark=triangle*,
    blue
] file[] {Adv_P1_Error.csv}; 

\addplot[
only marks,
mark=square,
    blue
] file[] {Adv_P2_Error.csv}; 
 
\addplot[
only marks,
mark=o,
    red
] file[] {Adv_P3_Error_80_alpha_beta_gamma.csv};

\addplot[
mark=diamond*,
    black,
    dashed
] file[] {Adv_P4_Error_80_alpha_beta_gamma.csv};

\addplot[
mark=o,
    black,
    solid
] file[] {WenoSofiaShock.csv};
  
\legend{
   P1 $\epst_1=\frac{1}{2}$, P2 $\epst_1=\frac{1}{6}$ $\epst_2=\frac{1}{25}$, P3 $\epst_1=\frac{1}{10}$ $\epst_2=\frac{1}{200}$ ,P4 $\epst_1=\frac{3}{40}$ $\epst_2=\frac{1}{500}$, WENO3}
\end{axis}
\end{tikzpicture}
\caption{Linear propagating discontinuity: (Top) Initial and final exact solution vs.  CG (Bottom) Error of CG as compared to third order WENO \cite{Sofia2011,Weno3_1996}.} 
\label{fig:adv_exactvsCG}
\end{figure}
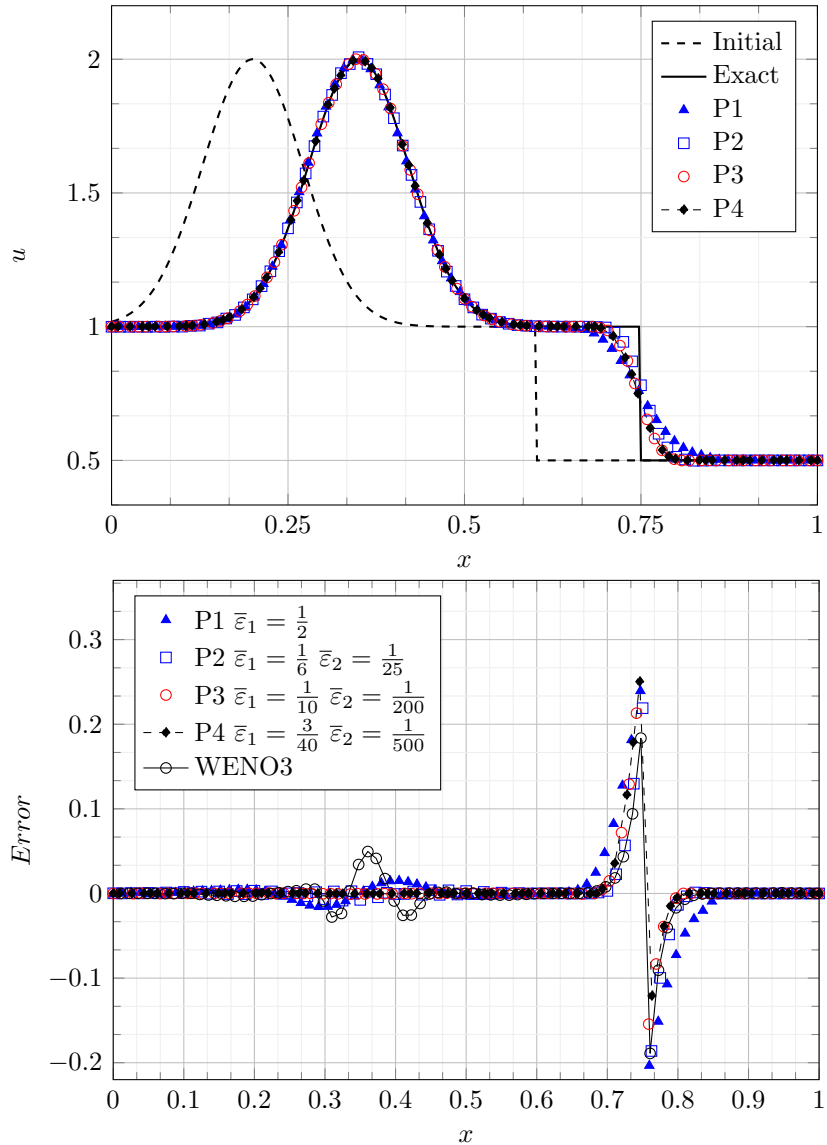

\begin{figure}[ht!]
\centering
\begin{tikzpicture}[font=\footnotesize]
 
\begin{axis}[
    xmin = 0.4, xmax = 1,
    ymin = -0.2, ymax = 0.25,
    xtick distance = 0.1,
    ytick distance = 0.1,
    legend pos = north west,
    legend cell align={left},
    grid = both,
    minor tick num = 2,
    major grid style = {lightgray},
    minor grid style = {lightgray!25},
    width = 0.8 \textwidth,
    height = 0.6 \textwidth,
    xlabel = {$x$},
    ylabel = {$Error$},]
 

\addplot[
     thick,
    red,
    solid
] file[] {Adv_P3_Error_80_alphadiv2.csv};
 
\addplot[
     thick,
     dashed,
blue
] file[] {Adv_P3_Error_80_alpha_beta_gamma.csv}; 
 
\addplot[
only marks,
mark=triangle*,
    red,
    solid
] file[] {Adv_P3_Error_80_alphax2.csv};

\legend{
   P3 $\epst_1=\frac{1}{20}$ $\epst_2=\frac{1}{200}$, P3 $\epst_1=\frac{1}{10}$ $\epst_2=\frac{1}{200}$,P3 $\epst_1=\frac{1}{5}$  $\epst_2=\frac{1}{200}$}
\end{axis}
\end{tikzpicture}

\begin{tikzpicture}[font=\footnotesize]
\begin{axis}[
    xmin = 0.4, xmax = 1,
    ymin = -0.233333, ymax = 0.23333,
    xtick distance = 0.1,
    ytick distance = 0.1,
    legend pos = north west,
    legend cell align={left},
    grid = both,
    minor tick num = 2,
    major grid style = {lightgray},
    minor grid style = {lightgray!25},
    width = 0.8 \textwidth,
    height = 0.6 \textwidth,
    xlabel = {$x$},
    ylabel = {$Error$},]
 

\addplot[
only marks,
mark=oplus,
    red,
] file[] {Adv_P3_Error_60_alpha_beta_gamma.csv};
 
\addplot[
only marks,
   mark=square,
    blue,
] file[] {Adv_P3_Error_80_alpha_beta_gamma.csv}; 
 
\addplot[
thick,
    blue,
    solid
] file[] {Adv_P3_Error_151_alpha_beta_gamma.csv};

\legend{
   P3 $60$ nodes, P3 $80$ nodes, P3 $151$ nodes}
\end{axis}
\end{tikzpicture}
 \caption{Linear propagating discontinuity: (Top) Error for CG P3 with varying amounts of first order dissipation. (Bottom) Error with $\epst_1=\frac{1}{10}$, $\epst_2=\frac{1}{200}$ for 3 mesh resolutions.}
\label{fig:adv_error_2}
\end{figure}
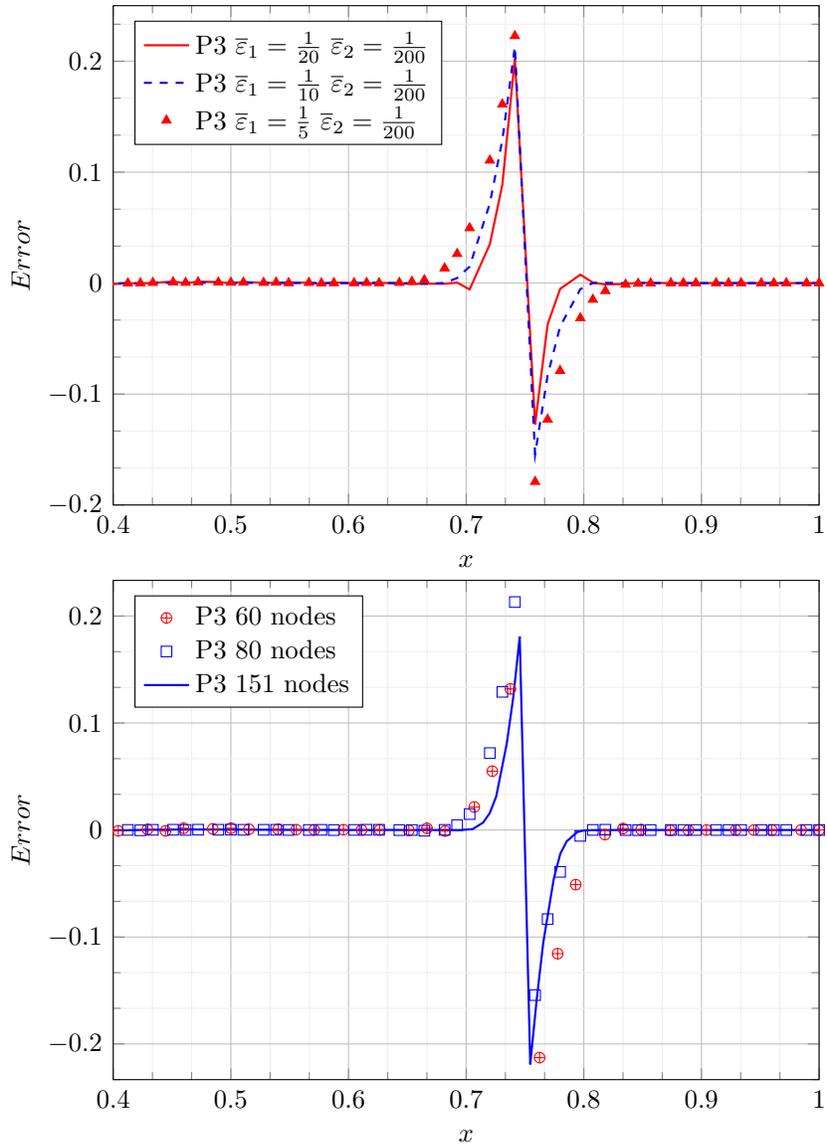

\begin{figure}[ht!]
\centering
\begin{tikzpicture}[font=\footnotesize]
 
\begin{axis}[
    xmin = 0.0, xmax = 0.6,
    ymin = -0.005, ymax = 0.0125,
    xtick distance = 0.1,
    ytick distance = 0.005,
    legend pos = north west,
    legend cell align={left},
    grid = both,
    minor tick num = 2,
    major grid style = {lightgray},
    minor grid style = {lightgray!25},
    width = 0.8 \textwidth,
    height = 0.6 \textwidth,
    xlabel = {$x$},
    ylabel = {$Error$},]
 

\addplot[
    thick,
    blue,
    solid
] file[] {Adv_P3_Error_80_alpha.csv}; 

\addplot[
only marks,
mark=diamond*,
    black
] file[] {Adv_P3_Error_80_alpha_beta.csv}; 
 
\addplot[
only marks,
mark=square,
red
] file[] {Adv_P3_Error_80_alpha_beta_gamma.csv};

\addplot[
mark=o,
    black
] file[] {WenoSofiaShock.csv};
  
\legend{
   P3 $\epst_1=\frac{1}{10}$, P3 $\epst_1=\frac{1}{10}$ $\epst_2=\frac{1}{200}$, P3 $\epst_1=\frac{1}{10}$ $\epst_2=\frac{1}{200}$ $\epst_3=\frac{1}{1000}$, WENO3}
\end{axis}
\end{tikzpicture}
\caption{Linear propagating discontinuity: Error of CG with different AD orders as compared to third order WENO \cite{Sofia2011,Weno3_1996}.} 
\label{fig:adv_error_3}
\end{figure}
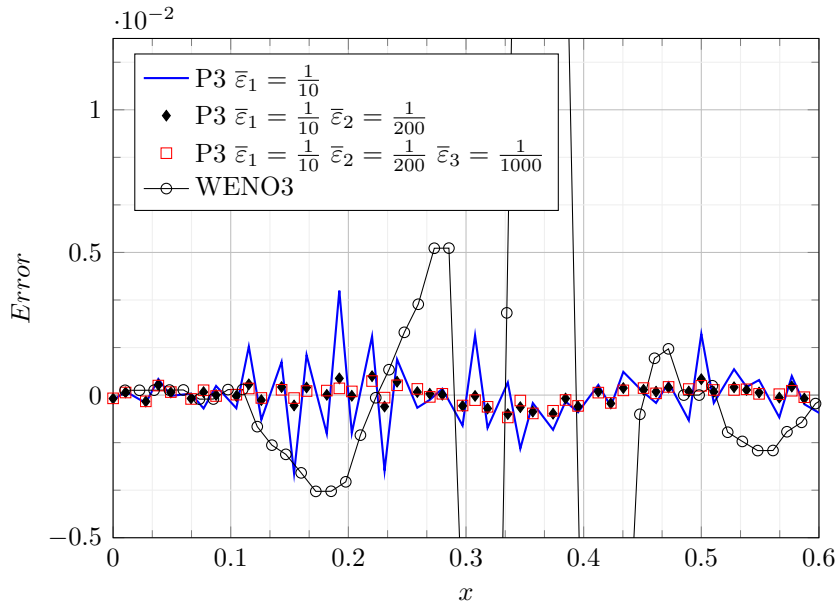

\subsection{Linear Smooth and Discontinuous Propagation}

We next consider the linear advection equation:
\begin{align}
 u_t + a u_x &= 0;~~~\Vol \in [0,1];~~~u(t,0)  =  u_0(x). \nn 
\end{align}

\noindent As per \cite{Sofia2011} the initial solution is a Gaussian pulse combined with a step function i.e.
\begin{equation}
u_0(x) = 
     \begin{cases}
       e^{-100(x-0.2)^2} & \quad\text{if }x \le 0.6\\
       0.5 &\quad\text{if } 0.6 <x \le 1.0 \\
     \end{cases}
\end{equation}

\noindent as shown in Figure~\ref{fig:adv_exactvsCG} (top). This case assesses the scheme's performance in modelling a linear problem in the presence of \textit{both} a smooth and discontinuous field.  We also seek to assess the performance of the dissipation term across a propagating discontinuity with a constant jump condition.

For temporal discretization, the Euler backward difference method with a sufficiently small time-step size is employed (set such that the temporal discretization error is vanishingly small compared to the spatial discretization error).  Following \cite{Sofia2011}, we employ a  $80$ node mesh and solve using CG P1 to P4.  The same $\epst_i$ was employed at nodes which lie with-in a distance of $0.1$ from the discontinuity, and set to zero elsewhere.  The resulting dissipation coefficients were found by trial-and-error to minimize spurious oscillations.
 
The solutions obtained with CG P1 to P4 are compared to 3rd order accurate WENO \cite{Weno3_1996} in Figure~\ref{fig:adv_exactvsCG}. Also shown is the corresponding error and AD coefficients that were used.  The P2, P3 and P4 solutions offer similar sharpness to WENO around the shock,  but outperform the latter over the pulse.  A new finding of this article is therefore the ability of the proposed CG formulation to retain high-order accuracy over the smooth field while at the same time achieving WENO like shock capturing over the discontinuity.

\begin{figure}[ht!]
\centering
\begin{tikzpicture}[font=\footnotesize]

\begin{axis}[
    xmin = 0.0, xmax = 1,
    ymin = -1.2, ymax = 1.2,
    xtick distance = 0.25,
    ytick distance = 0.5,
    legend pos = north east,
    legend cell align={left},
    legend style={font=\footnotesize},
    grid = both,
    minor tick num = 2,
    major grid style = {lightgray},
    minor grid style = {lightgray!25},
    width = 0.8\textwidth,
    height = 0.6\textwidth,
    xlabel = {$x$},
    ylabel = {$u$},]
 

\addplot[
    thick, dashed,
    black,
] file[] {SofiaBurgers_t0.001_P4_U_nodiss_81node.csv}; 

\addplot[
    thick, solid,
    black,
] file[] {SofiaBurgers_t0.1_V.csv};

\addplot[
only marks,
mark=square,
    blue
] file[] {SofiaBurgers_t0.1_P2_U_nodiss_81node.csv}; 

\addplot[
only marks,
mark=o,
    red
] file[] {SofiaBurgers_t0.1_P3_U_nodiss_79node.csv}; 

\addplot[
mark=diamond*,
    black,
    dashed
] file[] {SofiaBurgers_t0.1_P4_U_nodiss_81node.csv};

\legend{Initial,Exact,P2, P3, P4}

\end{axis}
\end{tikzpicture}

\begin{tikzpicture}[font=\footnotesize]
 
\begin{axis}[
    xmin = 0.0, xmax = 1,
    ymin = -0.015, ymax = 0.015,
    xtick distance = 0.25,
    ytick distance = 0.01,
    legend pos = north east,
    legend cell align={left},
    grid = both,
    minor tick num = 2,
    major grid style = {lightgray},
    minor grid style = {lightgray!25},
    width = 0.8 \textwidth,
    height = 0.6 \textwidth,
    xlabel = {$x$},
    ylabel = {$Error$},]
 

\addplot[
only marks,
mark=square,
    blue
] file[] {SofiaBurgers_t0.1_P2_E_nodiss_81node.csv}; 

\addplot[
only marks,
mark=o,
    red
] file[] {SofiaBurgers_t0.1_P3_E_nodiss_79node.csv}; 
 
\addplot[
    mark=diamond*,
    black,
    dashed
] file[] {SofiaBurgers_t0.1_P4_E_nodiss_81node.csv};

\addplot[
mark=o,
    black,
    solid
] file[] {Sofia_Burgers_error_0.1s_weno.csv};

\legend{
   P2 , P3,  P4, WENO3}
\end{axis}
\end{tikzpicture}
\caption{Non-linear problem: (Top) Initial and  exact solution at $t=0.1s$ for  CG schemes. (Bottom) Error of CG schemes as compared to third order WENO \cite{Sofia2011,Weno3_1996}.} 
\label{fig:nonlin_1}
\end{figure}
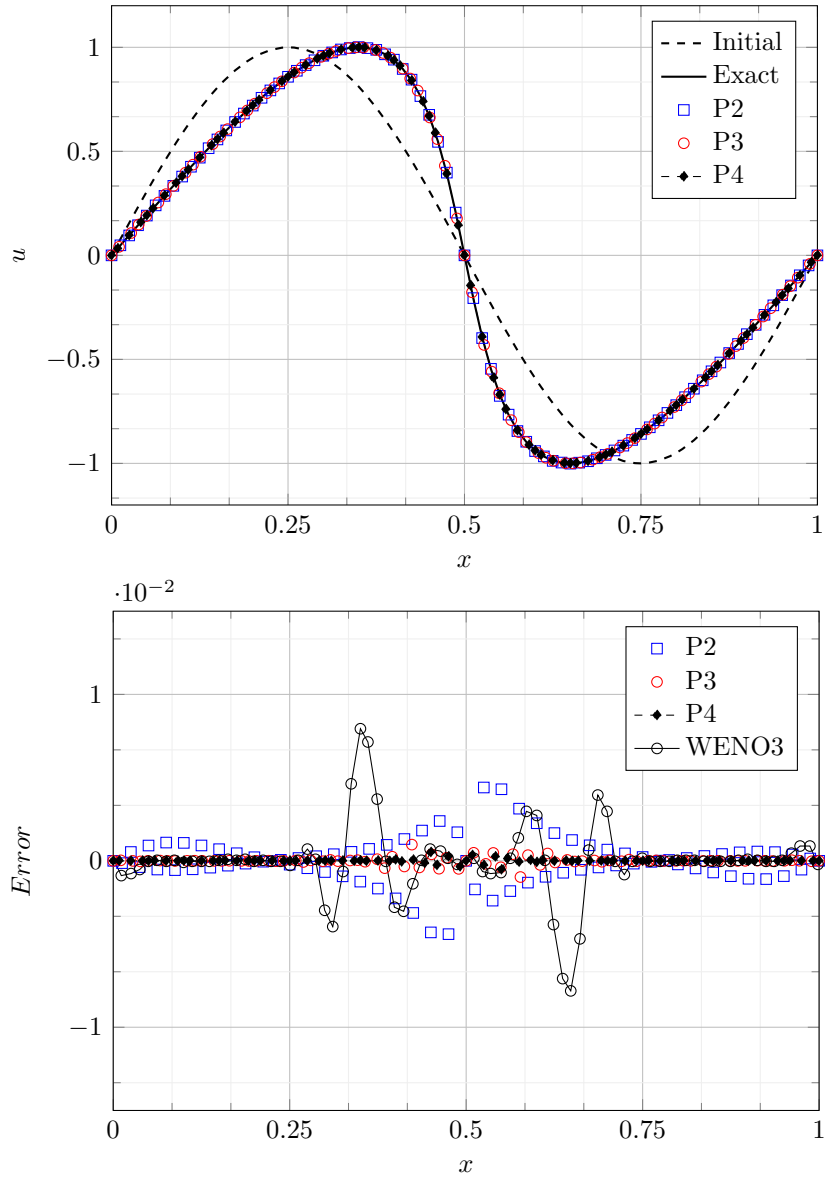 

The major contributor to ensuring sharpness in shock capturing was the value of $\epst_1$.  A too large or small value resulted in either under and over shoots or a less sharp shock.  This is illustrated for the P3 case in Figure~\ref{fig:adv_error_2}.  Also depicted in the figure (bottom) is the insensitivity of the chosen dissipation coefficients $\epst_1$ and $\epst_2$ to mesh resolution.  As shown, the dissipation effect is mesh spacing invariant as postulated in Section \ref{sec:paraCG}. This is evident from the oscillation free solutions for the three different meshes using the same dissipation coefficients. Figure~\ref{fig:adv_error_3} depicts the relative contribution of the higher order dissipation terms for P3. As shown the $\epst_2$ term has a significant effect upstream of the shock with the contribution of the $\epst_3$ term being significantly smaller. In all cases CG performs significantly better than WENO3.

\begin{figure}[ht!]
\centering
\begin{tikzpicture}[font=\footnotesize]

\begin{axis}[
    xmin = 0.0, xmax = 1,
    ymin = -1.2, ymax = 1.2,
    xtick distance = 0.25,
    ytick distance = 0.5,
    legend pos = north east,
    legend cell align={left},
    legend style={font=\footnotesize},
    grid = both,
    minor tick num = 2,
    major grid style = {lightgray},
    minor grid style = {lightgray!25},
    width = 0.8\textwidth,
    height = 0.6\textwidth,
    xlabel = {$x$},
    ylabel = {$u$},]
 

\addplot[
    thick, dashed,
    black,
] file[] {SofiaBurgers_t0.001_P4_U_nodiss_81node.csv}; 

\addplot[
    thick, solid,
    black,
] file[] {SofiaBurgers_t0.16_V.csv}; 

\addplot[
only marks,
mark=square,
    blue
] file[] {SofiaBurgers_t0.16_P2_U_nodiss_81node.csv}; 

\addplot[
only marks,
mark=o,
    red
] file[] {SofiaBurgers_t0.16_P3_U_nodiss_79node.csv}; 

\addplot[
   mark=diamond*,
    black,
    dashed
] file[] {SofiaBurgers_t0.16_P4_U_nodiss_81node.csv};

\legend{Initial,Exact,P2, P3, P4}

\end{axis}
\end{tikzpicture}

\begin{tikzpicture}[font=\footnotesize]
 
\begin{axis}[
    xmin = 0.0, xmax = 1,
    ymin = -0.15, ymax = 0.15,
    xtick distance = 0.25,
    ytick distance = 0.1,
    legend pos = north east,
    legend cell align={left},
    grid = both,
    minor tick num = 2,
    major grid style = {lightgray},
    minor grid style = {lightgray!25},
    width = 0.8 \textwidth,
    height = 0.6 \textwidth,
    xlabel = {$x$},
    ylabel = {$Error$},]
 

\addplot[
 only marks,
mark=square,
    blue
] file[] {SofiaBurgers_t0.16_P2_E_nodiss_81node.csv}; 

\addplot[
only marks,
mark=o,
    red
] file[] {SofiaBurgers_t0.16_P3_E_nodiss_79node.csv}; 
 
\addplot[
    mark=diamond*,
    black,
    dashed
] file[] {SofiaBurgers_t0.16_P4_E_nodiss_81node.csv};

\addplot[
    mark=o,
    black,
] file[] {Sofia_Burgers_error_0.16s_weno.csv};

\legend{
   P2 $\epst_1=\frac{1}{6}$ $\epst_2=\frac{1}{50}$, P3 $\epst_1=\frac{1}{8}$ $\epst_2=\frac{1}{500}$,P4 $\epst_2=\frac{1}{1000}$,WENO3}
\end{axis}
\end{tikzpicture}
\caption{Non-linear problem: (Top) Initial and  exact solution at $t=0.16s$ for  CG schemes. (Bottom) Error of CG schemes as compared to third order WENO \cite{Sofia2011,Weno3_1996}.} 
\label{fig:nonlin_2}
\end{figure}
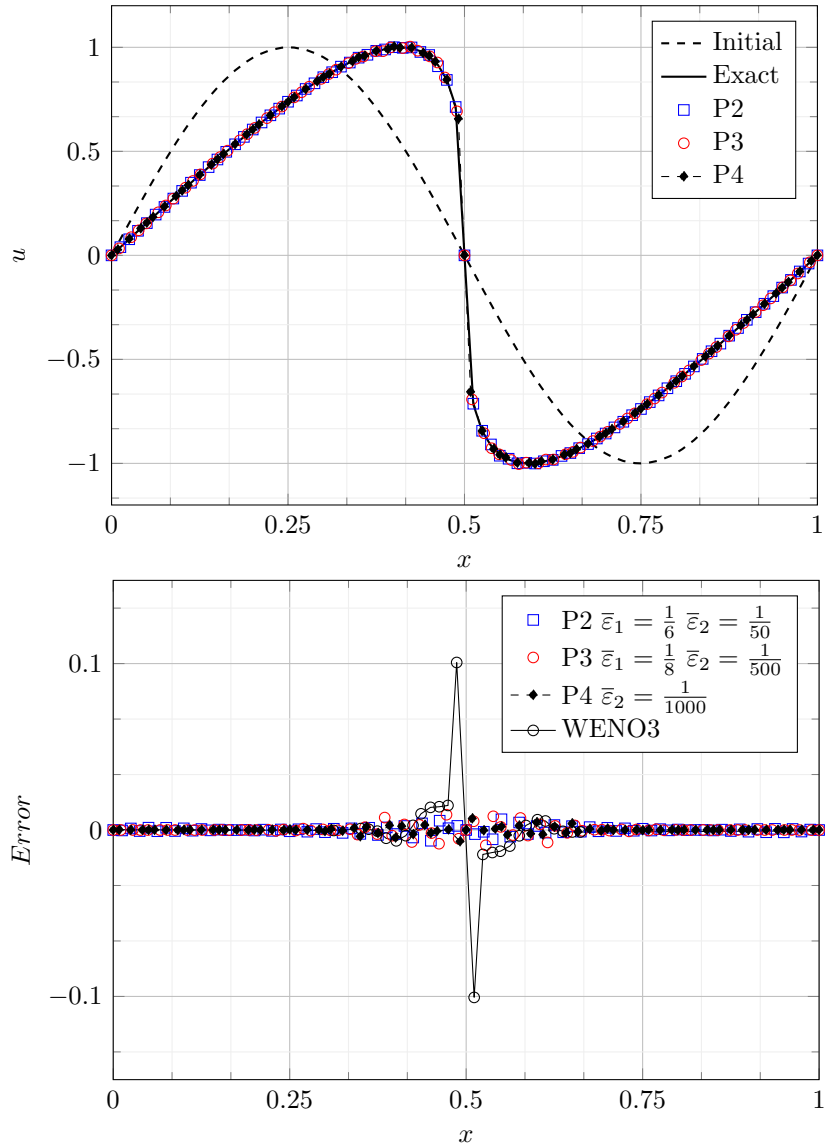  

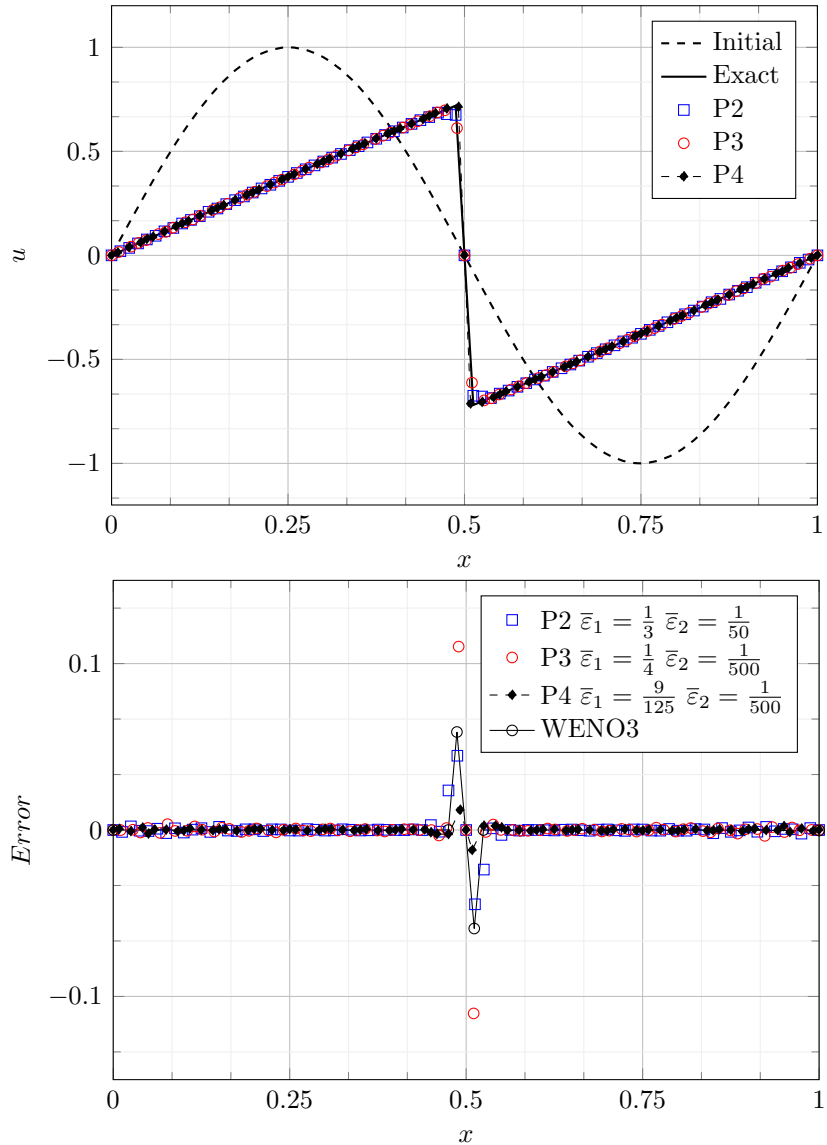
\begin{figure}[ht!]
\centering
\begin{tikzpicture}[font=\footnotesize]

\begin{axis}[
    xmin = 0.0, xmax = 1,
    ymin = -1.2, ymax = 1.2,
    xtick distance = 0.25,
    ytick distance = 0.5,
    legend pos = north east,
    legend cell align={left},
    legend style={font=\footnotesize},
    grid = both,
    minor tick num = 2,
    major grid style = {lightgray},
    minor grid style = {lightgray!25},
    width = 0.8\textwidth,
    height = 0.6\textwidth,
    xlabel = {$x$},
    ylabel = {$u$},]
 

\addplot[
    thick, dashed,
    black,
] file[] {SofiaBurgers_t0.001_P4_U_nodiss_81node.csv}; 

\addplot[
    thick, solid,
    black,
] file[] {SofiaBurgers_t0.5_V.csv}; 

\addplot[
only marks,
mark=square,
    blue
] file[] {SofiaBurgers_t0.5_P2_U_diss_81node.csv}; 

\addplot[
only marks,
mark=o,
    red
] file[] {SofiaBurgers_t0.5_P3_U_diss_79node.csv}; 

\addplot[
    mark=diamond*,
    black,
    dashed
] file[] {SofiaBurgers_t0.5_P4_U_diss_81node.csv};

\legend{Initial,Exact,P2, P3, P4}

\end{axis}
\end{tikzpicture}

\begin{tikzpicture}[font=\footnotesize]
 
\begin{axis}[
    xmin = 0.0, xmax = 1,
    ymin = -0.15, ymax = 0.15,
    xtick distance = 0.25,
    ytick distance = 0.1,
    legend pos = north east,
    legend cell align={left},
    grid = both,
    minor tick num = 2,
    major grid style = {lightgray},
    minor grid style = {lightgray!25},
    width = 0.8 \textwidth,
    height = 0.6 \textwidth,
    xlabel = {$x$},
    ylabel = {$Error$},]
 

\addplot[
only marks,
mark=square,
    blue
] file[] {SofiaBurgers_t0.5_P2_E_diss_81node.csv}; 

\addplot[
only marks,
mark=o,
    red
] file[] {SofiaBurgers_t0.5_P3_E_diss_79node.csv}; 
 
\addplot[
    mark=diamond*,
    black,
    dashed
] file[] {SofiaBurgers_t0.5_P4_E_diss_81node.csv};

\addplot[
    mark=o,
    black,
] file[] {Sofia_Burgers_error_0.5s_weno.csv};

\legend{
   P2 $\epst_1=\frac{1}{3}$ $\epst_2=\frac{1}{50}$, P3 $\epst_1=\frac{1}{4}$ $\epst_2=\frac{1}{500}$,P4 $\epst_1=\frac{9}{125}$ $\epst_2=\frac{1}{500}$,WENO3}
\end{axis}
\end{tikzpicture}
\caption{Non-linear problem: (Top) Initial and  exact solution at $t=0.5s$ for  CG schemes. (Bottom) Error of CG schemes as compared to third order WENO \cite{Sofia2011,Weno3_1996}.} 
\label{fig:nonlin_3}
\end{figure} 

\subsection{Non-Linear Smooth and Discontinuous Propagation}

We finally consider the Burger's equation:
\begin{align}\label{eq:non_lin}
 u_t + u u_x &= 0;~~~\Vol \in [0,1]; \nn 
\end{align}

\noindent where as per \cite{Sofia2011} the initial solution is taken as $u(x,0)=sin(2 \pi x)$.  The analytical solution is found from the roots of $\left[ -u + sin(2 \pi (-u~t + x)) \right]=0$.  The initial condition is shown in Figure~\ref{fig:nonlin_1}.

\begin{figure}[ht!]
\centering
\begin{tikzpicture}[font=\footnotesize]
 
\begin{axis}[
    xmin = 0.0, xmax = 1,
    ymin = -0.1, ymax = 0.15,
    xtick distance = 0.25,
    ytick distance = 0.1,
    legend pos = north east,
    legend cell align={left},
    grid = both,
    minor tick num = 2,
    major grid style = {lightgray},
    minor grid style = {lightgray!25},
    width = 0.8 \textwidth,
    height = 0.6 \textwidth,
    xlabel = {$x$},
    ylabel = {$Error$},]
 

\addplot[
only marks,
mark=square,
    blue
] file[] {SofiaBurgers_t0.5_P4_E_diss_77node.csv}; 
 
\addplot[
    mark=diamond*,
dashed,
] file[] {SofiaBurgers_t0.5_P4_E_diss_81node.csv};

\addplot[
mark=o,
    black,
    solid
] file[] {Sofia_Burgers_error_0.5s_weno.csv};

\legend{
   P4: 19 elements (77 nodes), P4: 20 elements (81 nodes), WENO3: 80 nodes}
\end{axis}
\end{tikzpicture}
\caption{Non-linear problem: Error of P4 using odd and even number of elements as compared to third order WENO \cite{Sofia2011,Weno3_1996}.} 
\label{fig:nonlin_4}
\end{figure}
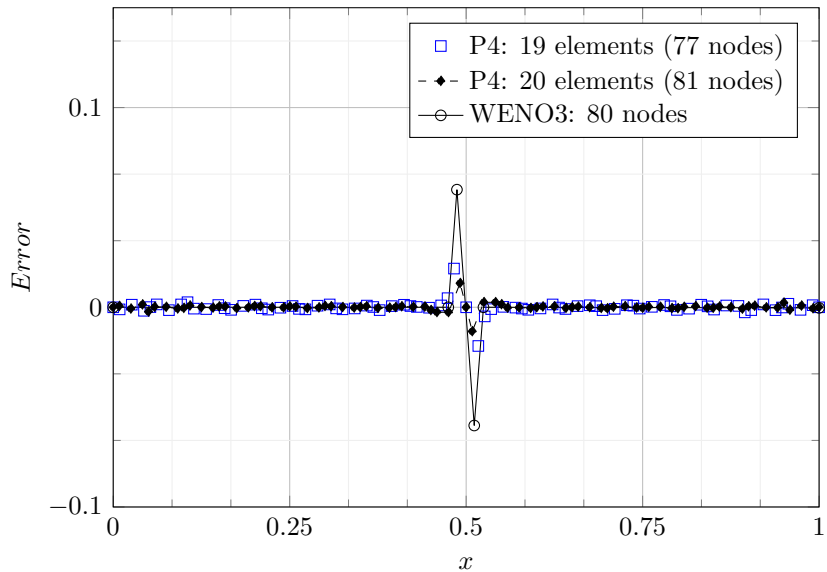 

The governing equation is discretized in the following split form \cite{JanNonLin2013} for stability
\begin{equation} \label{eq:non_lin_d}
\Pp \U_t + \frac{1}{3} \left[ \U \Qx \U + \Qx \U^2 \right] + |\U|  \sum_e \Di^e \U = 0.
\end{equation}

Applying dissipation to the non-linear problem is more challenging as the field $u$ evolves from continuous to discontinuous over time as per Figures~\ref{fig:nonlin_1} to \ref{fig:nonlin_3}.  Dissipation is therefore only activated over the discontinuity ($ 0.4<x<0.6$) for $t > 0.15$ (once the discontinuity starts to appear).  We again use spatially constant dissipation coefficients $\epst_1,\epst_2$ over this range and set $\epst_i=0$  elsewhere.

We consider next the accuracy of the CG formulation in Figures~\ref{fig:nonlin_1} and \ref{fig:nonlin_2}.  As shown CG significantly outperforms the third order WENO while the P3 and P4 cases are comparable to the $4^{th}$ order accurate FD results in \cite{Sofia2011}.  Note that the small spurious oscillations are due insufficient mesh resolution (as discussed in \cite{Frenander2018}) and not an instability. 

Finally, we consider $t=0.5$ at which point the discontinuity is fully developed (Figure~\ref{fig:nonlin_3}).  Larger amounts of dissipation than previously is now required.  As shown, $p=2$ offers accuracy over the shock which is comparable to WENO3 with $p=3$ being less accurate. However, $p=4$ significantly outperforms WENO. To demonstrate that this is invariant with respect to the position of the discontinuity relative to an element, the solution is compared to WENO on an odd and even number of elements in Figure~\ref{fig:nonlin_4}. This demonstrates that CG with $p=4$ and an added suitable AD can be used as an effective shock capturing scheme.

\section{Conclusions}
This article presented a high-order continuous Galerkin based FE scheme for the purpose of solving the advection-diffusion and wave propagation equations in the presence of both smooth and discontinuous fields.  Keys to developing a provably stable scheme were the SBP-SAT formulation as well as the development of provably stable Galerkin weighted artificial dissipation operators that remove spurious oscillations over discontinuous fields.  The new artificial dissipation operators may also be used without modification in DG schemes. The developed CG formulation was proven stable and exemplified via a series of numerical experiments which employed Lagrange polynomials  of order $p=1$ to $p=4$.  For the smooth advection-diffusion equation super convergence of order $(p+2)$ was achieved.  For both linear and non-linear transient advection equations, WENO like accuracy was demonstrated over shocks while high-order accuracy was retained elsewhere. Due these findings the developed SBP-SAT CG scheme is deemed of interest for further investigation.

\section{Acknowledgements}
This work was financially supported by the European Union’s Horizon 2020 research and innovation programme under grant agreement No 81504 (the SLOshing Wing Dynamics (SLOWD) project) as well as the National Research Foundation (NRF) of South Africa grant 89916. The opinions, findings, and conclusions or recommendations expressed are those of the authors alone, and the NRF accepts no liability whatsoever in this regard. Jan Nordström was funded by Vetenskapsrådet Sweden grant 2021-05484 VR and the University of Johannesburg, South Africa.

\bibliographystyle{model1-num-names}
\bibliography{Bibliography}

\begin{thebibliography}{47}
\expandafter\ifx\csname natexlab\endcsname\relax\def\natexlab#1{#1}\fi
\providecommand{\bibinfo}[2]{#2}
\ifx\xfnm\relax \def\xfnm[#1]{\unskip,\space#1}\fi
\bibitem[{Kreiss and Oliger(1972)}]{Kreiss1972}
\bibinfo{author}{H.~Kreiss}, \bibinfo{author}{J.~Oliger},
\newblock \bibinfo{title}{Comparison of accurate methods for the integration of
  hyperbolic equations},
\newblock volume \bibinfo{volume}{TellusXXIV(3) (1972)}.
\bibitem[{Fernandez et~al.(2017)Fernandez, Nguyen, and Peraire}]{Peraire2017}
\bibinfo{author}{P.~Fernandez}, \bibinfo{author}{N.~C. Nguyen},
  \bibinfo{author}{J.~Peraire},
\newblock \bibinfo{title}{The hybridized discontinuous \uppercase{G}alerkin
  method for implicit large-eddy simulation of transitional turbulent flows},
\newblock \bibinfo{journal}{Journal of Computational Physics}
  \bibinfo{volume}{336} (\bibinfo{year}{2017}) \bibinfo{pages}{308--329}.
\bibitem[{Loppi et~al.(2018)Loppi, Witherden, Jameson, and
  Vincent}]{Jameson2018}
\bibinfo{author}{N.~A. Loppi}, \bibinfo{author}{F.~D. Witherden},
  \bibinfo{author}{A.~Jameson}, \bibinfo{author}{P.~E. Vincent},
\newblock \bibinfo{title}{A high-order cross-platform incompressible
  navier-stokes solver via artificial compressibility with application to a
  turbulent jet},
\newblock \bibinfo{journal}{Computer Physics CommunicationsS}
  \bibinfo{volume}{233} (\bibinfo{year}{2018}) \bibinfo{pages}{193--205}.
\bibitem[{Changfoot et~al.(2019)Changfoot, Malan, and
  Nordstr\"om}]{Changfoot19}
\bibinfo{author}{D.~M. Changfoot}, \bibinfo{author}{A.~G. Malan},
  \bibinfo{author}{J.~Nordstr\"om},
\newblock \bibinfo{title}{{Hybrid Computational-Fluid-Dynamics Platform to
  Investigate Aircraft Trailing Vortices}},
\newblock \bibinfo{journal}{{Journal of Aircraft}} \bibinfo{volume}{{56}}
  (\bibinfo{year}{{2019}}) \bibinfo{pages}{{344--355}}.
\bibitem[{Ilangakoon et~al.(2020)Ilangakoon, Malan, and Jones}]{Ilangakoon2020}
\bibinfo{author}{N.~A. Ilangakoon}, \bibinfo{author}{A.~G. Malan},
  \bibinfo{author}{B.~W.~S. Jones},
\newblock \bibinfo{title}{A higher-order accurate surface tension modelling
  volume-of-fluid scheme for 2d curvilinear meshes},
\newblock \bibinfo{journal}{Journal of Computational Physics}
  \bibinfo{volume}{420} (\bibinfo{year}{2020}).
\bibitem[{Vila-Perez et~al.(2021)Vila-Perez, Giacomini, Sevilla, and
  Huerta}]{Sevilla2021}
\bibinfo{author}{J.~Vila-Perez}, \bibinfo{author}{M.~Giacomini},
  \bibinfo{author}{R.~Sevilla}, \bibinfo{author}{A.~Huerta},
\newblock \bibinfo{title}{Hybridisable discontinuous \uppercase{G}alerkin
  formulation of compressible flows},
\newblock \bibinfo{journal}{Archives of Computational Methods in Engineering}
  \bibinfo{volume}{28} (\bibinfo{year}{2021}) \bibinfo{pages}{753--784}.
\bibitem[{Ilangakoon and Malan(2022)}]{Ilangakoon2022}
\bibinfo{author}{N.~A. Ilangakoon}, \bibinfo{author}{A.~G. Malan},
\newblock \bibinfo{title}{A higher-order accurate vof interface curvature
  computation scheme for 3d non-orthogonal structured meshes},
\newblock \bibinfo{journal}{Computers \& Fluids} \bibinfo{volume}{245}
  (\bibinfo{year}{2022}).
\bibitem[{Nordstr\"om(2023)}]{Nordstrom2023}
\bibinfo{author}{J.~Nordstr\"om},
\newblock \bibinfo{title}{Nonlinear boundary conditions for initial boundary
  value problems with applications in computational fluid dynamics},
\newblock \bibinfo{journal}{arXiv:2306.01297}  (\bibinfo{year}{2023}).
\bibitem[{Huerta et~al.(2013)Huerta, Angeloski, Roca, and Peraire}]{Huerte2013}
\bibinfo{author}{A.~Huerta}, \bibinfo{author}{A.~Angeloski},
  \bibinfo{author}{X.~Roca}, \bibinfo{author}{J.~Peraire},
\newblock \bibinfo{title}{Efficiency of high-order elements for continuous and
  discontinuous \uppercase{G}alerkin methods},
\newblock \bibinfo{journal}{International Journal for Numerical Methods in
  Engineering} \bibinfo{volume}{96} (\bibinfo{year}{2013})
  \bibinfo{pages}{529--560}.
\bibitem[{Kreiss and Scherer(1974)}]{Kreiss1974}
\bibinfo{author}{H.~Kreiss}, \bibinfo{author}{G.~Scherer},
\newblock \bibinfo{title}{Finite element and finite difference methods for
  hyperbolic partial differential equations, mathematical aspects of finite
  elements in partial differential equations paper},
\newblock in: \bibinfo{booktitle}{Proceedings of a Symposium Conducted by the
  Mathematics Research Center, the University of Wisconsin–Madison (1974)}.
\bibitem[{Kreiss and Scherer(1977)}]{Kreiss1977}
\bibinfo{author}{H.~Kreiss}, \bibinfo{author}{G.~Scherer},
\newblock \bibinfo{title}{On the existence of energy estimates for difference
  approximations for hyperbolic systems},
\newblock in: \bibinfo{booktitle}{Technicalreport,Dept.ofScientific Computing,
  Uppsala University (1977)}.
\bibitem[{Carpenter et~al.(1994)Carpenter, Gottlieb, and
  Abarbanel}]{Carpenter1994}
\bibinfo{author}{M.~Carpenter}, \bibinfo{author}{D.~Gottlieb},
  \bibinfo{author}{S.~Abarbanel},
\newblock \bibinfo{title}{{Time-Stable Boundary-Conditions for
  Finite-Difference Schemes Solving Hyperbolic Systems - Methodology and
  Application to High-Order Compact Schemes}},
\newblock \bibinfo{journal}{{Journal of Computational Physics}}
  \bibinfo{volume}{{111}} (\bibinfo{year}{{1994}}) \bibinfo{pages}{{220--236}}.
\bibitem[{Carpenter et~al.(1999)Carpenter, Nordstr\"om, and
  Gottlieb}]{Carpenter1991}
\bibinfo{author}{M.~Carpenter}, \bibinfo{author}{J.~Nordstr\"om},
  \bibinfo{author}{D.~Gottlieb},
\newblock \bibinfo{title}{A stable and conservative interface treatment of
  arbitrary spatial accuracy},
\newblock \bibinfo{journal}{Journal of Computational Physics}
  \bibinfo{volume}{148} (\bibinfo{year}{1999}) \bibinfo{pages}{341--365}.
\bibitem[{Nordstr\"om and Carpenter(1999)}]{Nordstrom1999}
\bibinfo{author}{J.~Nordstr\"om}, \bibinfo{author}{M.~Carpenter},
\newblock \bibinfo{title}{Boundary and interface conditions for high-order
  finite-difference methods applied to the \uppercase{E}uler and
  \uppercase{N}avier-stokes equations},
\newblock \bibinfo{journal}{Journal of Computational Physics}
  \bibinfo{volume}{148} (\bibinfo{year}{1999}) \bibinfo{pages}{621--645}.
\bibitem[{Nordstr\"om and Carpenter(2001)}]{Nordstrom2001}
\bibinfo{author}{J.~Nordstr\"om}, \bibinfo{author}{M.~Carpenter},
\newblock \bibinfo{title}{High-order finite difference methods,
  multidimensional linear problems, and curvilinear coordinates},
\newblock \bibinfo{journal}{Journal of Computational Physics}
  \bibinfo{volume}{173} (\bibinfo{year}{2001}) \bibinfo{pages}{149--174}.
\bibitem[{Sv\"ard and Nordstr\"om(2014)}]{Svard2014}
\bibinfo{author}{M.~Sv\"ard}, \bibinfo{author}{J.~Nordstr\"om},
\newblock \bibinfo{title}{Review of summation-by-parts schemes for
  initial-boundary-value problems},
\newblock \bibinfo{journal}{Journal of Computational Physics}
  \bibinfo{volume}{268} (\bibinfo{year}{2014}) \bibinfo{pages}{17--38}.
\bibitem[{Fernandez et~al.(2014)Fernandez, Hicken, and Zingg}]{Fernandez2020}
\bibinfo{author}{D.~C. D.~R. Fernandez}, \bibinfo{author}{J.~E. Hicken},
  \bibinfo{author}{D.~W. Zingg},
\newblock \bibinfo{title}{Review of summation-by-parts operators with
  simultaneous approximation terms for the numerical solution of partial
  differential equations},
\newblock \bibinfo{journal}{Computers \& Fluids} \bibinfo{volume}{95}
  (\bibinfo{year}{2014}) \bibinfo{pages}{171--196}.
\bibitem[{Abgrall et~al.(2020)Abgrall, Nordstr\"om, \"Offner, and
  Tokareva}]{2020JanCG1}
\bibinfo{author}{R.~Abgrall}, \bibinfo{author}{J.~Nordstr\"om},
  \bibinfo{author}{P.~\"Offner}, \bibinfo{author}{S.~Tokareva},
\newblock \bibinfo{title}{Analysis of the sbp-sat stabilization for finite
  element methods part i: Linear problems},
\newblock \bibinfo{journal}{Journal of Scientific Computing}
  \bibinfo{volume}{85} (\bibinfo{year}{2020}).
\bibitem[{Abgrall et~al.(2023)Abgrall, Nordstr\"om, \"Offner, and
  Tokareva}]{2023JanCG2}
\bibinfo{author}{R.~Abgrall}, \bibinfo{author}{J.~Nordstr\"om},
  \bibinfo{author}{P.~\"Offner}, \bibinfo{author}{S.~Tokareva},
\newblock \bibinfo{title}{Analysis of the sbp-sat stabilization for finite
  element methods part ii: Entropy stability},
\newblock \bibinfo{journal}{Communications on Applied Mathematics and
  Computation} \bibinfo{volume}{5} (\bibinfo{year}{2023})
  \bibinfo{pages}{573–595}.
\bibitem[{Jameson(1995)}]{Jameson1995}
\bibinfo{author}{A.~Jameson},
\newblock \bibinfo{title}{Analysis and design of numerical schemes for gas
  dynamics},
\newblock \bibinfo{journal}{International Journal of Computational Fluid
  Dynamics} \bibinfo{volume}{5} (\bibinfo{year}{1995}).
\bibitem[{Sv\"ard et~al.(2006)Sv\"ard, Gong, and Nordstr\"om}]{Nordstrom2006_2}
\bibinfo{author}{M.~Sv\"ard}, \bibinfo{author}{J.~Gong},
  \bibinfo{author}{J.~Nordstr\"om},
\newblock \bibinfo{title}{Stable artificial dissipation operators for finite
  volume schemes on unstructured grids},
\newblock \bibinfo{journal}{Applied Numerical Mathematics} \bibinfo{volume}{56}
  (\bibinfo{year}{2006}) \bibinfo{pages}{1481--1490}.
\bibitem[{Pattinson et~al.(2007)Pattinson, Malan, and Meyer}]{Pattinson2007}
\bibinfo{author}{J.~Pattinson}, \bibinfo{author}{A.~G. Malan},
  \bibinfo{author}{J.~P. Meyer},
\newblock \bibinfo{title}{A cut-cell non-conforming cartesian mesh method for
  compressible and incompressible flow},
\newblock \bibinfo{journal}{International Journal for Numerical Methods in
  Engineering} \bibinfo{volume}{72} (\bibinfo{year}{2007})
  \bibinfo{pages}{1332--1354}.
\bibitem[{Mattsson et~al.(2004)Mattsson, Sv\"ard, and
  Nordstr\"om}]{Mattsson2003}
\bibinfo{author}{K.~Mattsson}, \bibinfo{author}{M.~Sv\"ard},
  \bibinfo{author}{J.~Nordstr\"om},
\newblock \bibinfo{title}{Stable and accurate artificial dissipation},
\newblock \bibinfo{journal}{Journal of Scientific Computing}
  \bibinfo{volume}{21} (\bibinfo{year}{2004}) \bibinfo{pages}{57--79}.
\bibitem[{Lundquist and Nordstr\"om(2020)}]{Lundquist2020}
\bibinfo{author}{T.~Lundquist}, \bibinfo{author}{J.~Nordstr\"om},
\newblock \bibinfo{title}{Stable and accurate filtering procedures},
\newblock \bibinfo{journal}{Journal of Scientific Computing}
  \bibinfo{volume}{82} (\bibinfo{year}{2020}).
\bibitem[{Diener et~al.(2007)Diener, Dorband, Schnetter, and
  Tiglio}]{Diener2007}
\bibinfo{author}{P.~Diener}, \bibinfo{author}{E.~N. Dorband},
  \bibinfo{author}{E.~Schnetter}, \bibinfo{author}{M.~Tiglio},
\newblock \bibinfo{title}{Optimized high-order derivative and dissipation
  operators \uppercase{S}ummation-\uppercase{B}y-\uppercase{P}arts, and
  applications in three-dimensional multi-block evolutions},
\newblock \bibinfo{journal}{JOURNAL OF SCIENTIFIC COMPUTING}
  \bibinfo{volume}{32} (\bibinfo{year}{2007}) \bibinfo{pages}{109--145}.
\bibitem[{Sv\"ard and Mishra(2009)}]{Svard2009}
\bibinfo{author}{M.~Sv\"ard}, \bibinfo{author}{S.~Mishra},
\newblock \bibinfo{title}{Shock capturing artificial dissipation for high-order
  finite difference schemes},
\newblock \bibinfo{journal}{JOURNAL OF SCIENTIFIC COMPUTING}
  \bibinfo{volume}{39} (\bibinfo{year}{2009}) \bibinfo{pages}{454--484}.
\bibitem[{Penner and Zingg(2018)}]{Zingg2018}
\bibinfo{author}{D.~Penner}, \bibinfo{author}{D.~Zingg},
\newblock \bibinfo{title}{High-order artiﬁcial dissipation operators
  possessing the \uppercase{S}ummation-\uppercase{B}y-\uppercase{P}arts
  property},
\newblock in: \bibinfo{booktitle}{2018 Fluid Dynamics Conference, AIAA (2018)}.
\bibitem[{Ranocha et~al.(2018)Ranocha, Glaubitz, Oeffner, and
  Sonar}]{Ranocha2018}
\bibinfo{author}{H.~Ranocha}, \bibinfo{author}{J.~Glaubitz},
  \bibinfo{author}{P.~Oeffner}, \bibinfo{author}{T.~Sonar},
\newblock \bibinfo{title}{Stability of artificial dissipation and modal
  filtering for flux reconstruction schemes using summation-by-parts
  operators},
\newblock \bibinfo{journal}{APPLIED NUMERICAL MATHEMATICS}
  \bibinfo{volume}{128} (\bibinfo{year}{2018}) \bibinfo{pages}{1--23}.
\bibitem[{Tadmor(1989)}]{Tadmor1989}
\bibinfo{author}{E.~Tadmor},
\newblock \bibinfo{title}{Convergence of spectral methods for nonlinear
  conservation laws},
\newblock \bibinfo{journal}{SIAM J. Numer. Anal.}  (\bibinfo{year}{1989}).
\bibitem[{Maday et~al.(1993)Maday, Kaber, and Tadmor}]{Maday1993}
\bibinfo{author}{Y.~Maday}, \bibinfo{author}{S.~O. Kaber},
  \bibinfo{author}{E.~Tadmor},
\newblock \bibinfo{title}{Legendre pseudospectral viscosity method for
  nonlinear conservation laws},
\newblock \bibinfo{journal}{SIAM J. Numer. Anal.}  (\bibinfo{year}{1993}).
\bibitem[{Nguyen et~al.(2013)Nguyen, Vila-Pérez, and Peraire}]{Peraire2023}
\bibinfo{author}{N.~Nguyen}, \bibinfo{author}{J.~Vila-Pérez},
  \bibinfo{author}{J.~Peraire},
\newblock \bibinfo{title}{Discretely conservative finite-difference
  formulations for nonlinear conservation laws in split form: Theory and
  boundary conditions},
\newblock \bibinfo{journal}{Journal of Computational Physics}
  \bibinfo{volume}{234} (\bibinfo{year}{2013}) \bibinfo{pages}{353--375}.
\bibitem[{Hughes et~al.(1986)Hughes, Mallet, and Akira}]{Hughes1986}
\bibinfo{author}{T.~J. Hughes}, \bibinfo{author}{M.~Mallet},
  \bibinfo{author}{A.~Akira},
\newblock \bibinfo{title}{A new ﬁnite element formulation for computational
  ﬂuid dynamics: Ii. beyond \uppercase{SUPG}},
\newblock \bibinfo{journal}{Computer Methods in Applied Mechanics and
  Engineering} \bibinfo{volume}{54} (\bibinfo{year}{1986}).
\bibitem[{Manzari et~al.(1998)Manzari, Hassan, Morgan, and
  Weatherill}]{Manzari1998}
\bibinfo{author}{M.~Manzari}, \bibinfo{author}{O.~Hassan},
  \bibinfo{author}{K.~Morgan}, \bibinfo{author}{N.~Weatherill},
\newblock \bibinfo{title}{Turbulent flow computations on 3d unstructured
  grids},
\newblock \bibinfo{journal}{Finite Elements in Analysis and Design}
  \bibinfo{volume}{30} (\bibinfo{year}{1998}) \bibinfo{pages}{353--363}.
\bibitem[{Nordstr\"om(2017)}]{Nordstrom2017}
\bibinfo{author}{J.~Nordstr\"om},
\newblock \bibinfo{title}{A roadmap to well posed and stable problems in
  computational physics},
\newblock \bibinfo{journal}{Journal of Scientific Computing}
  \bibinfo{volume}{71} (\bibinfo{year}{2017}) \bibinfo{pages}{365--385}.
\bibitem[{Hesthaven and Gottlieb(1996)}]{Hesthaven1996}
\bibinfo{author}{J.~Hesthaven}, \bibinfo{author}{D.~Gottlieb},
\newblock \bibinfo{title}{A stable penalty method for the compressible
  \uppercase{N}avier-\uppercase{S}tokes equations .1. open boundary
  conditions},
\newblock \bibinfo{journal}{SIAM JOURNAL ON SCIENTIFIC COMPUTING}
  \bibinfo{volume}{17} (\bibinfo{year}{1996}) \bibinfo{pages}{579--612}.
\bibitem[{Gassner(2013)}]{Gassner2013}
\bibinfo{author}{G.~J. Gassner},
\newblock \bibinfo{title}{A skew-symmetric discontinuous \uppercase{G}alerkin
  spectral element discretization and its relation to \uppercase{SBP-SAT}
  finite difference methods},
\newblock \bibinfo{journal}{SIAM JOURNAL ON SCIENTIFIC COMPUTING}
  \bibinfo{volume}{35} (\bibinfo{year}{2013}) \bibinfo{pages}{A1233--A1253}.
\bibitem[{Kopriva et~al.(2021)Kopriva, Gassner, and Nordstr\"om}]{Kopriva2021}
\bibinfo{author}{D.~A. Kopriva}, \bibinfo{author}{G.~J. Gassner},
  \bibinfo{author}{J.~Nordstr\"om},
\newblock \bibinfo{title}{Stability of discontinuous \uppercase{G}alerkin
  spectral element schemes for wave propagation when the coefficient matrices
  have jumps},
\newblock \bibinfo{journal}{JOURNAL OF SCIENTIFIC COMPUTING}
  \bibinfo{volume}{88} (\bibinfo{year}{2021}).
\bibitem[{Arnold et~al.(2002)Arnold, Brezzi, Cockburn, and Marini}]{Brezzi01}
\bibinfo{author}{D.~Arnold}, \bibinfo{author}{F.~Brezzi},
  \bibinfo{author}{B.~Cockburn}, \bibinfo{author}{L.~Marini},
\newblock \bibinfo{title}{{Unified analysis of discontinuous
  \uppercase{G}alerkin methods for elliptic problems}},
\newblock \bibinfo{journal}{{SIAM Journal on Numerical Analysis}}
  \bibinfo{volume}{{39}} (\bibinfo{year}{{2002}})
  \bibinfo{pages}{{1749--1779}}.
\bibitem[{Nordstr\"om and Sv\"ard(2005)}]{Nordstrom05}
\bibinfo{author}{J.~Nordstr\"om}, \bibinfo{author}{M.~Sv\"ard},
\newblock \bibinfo{title}{{Well-posed boundary conditions for the Navier-Stokes
  equations}},
\newblock \bibinfo{journal}{{SIAM Journal on Numerical Analysis}}
  \bibinfo{volume}{{43}} (\bibinfo{year}{{2005}})
  \bibinfo{pages}{{1231--1255}}.
\bibitem[{Nordstr\"om(2006)}]{Nordstrom2006}
\bibinfo{author}{J.~Nordstr\"om},
\newblock \bibinfo{title}{Conservative finite difference formulations, variable
  coefficients, energy estimates and artificial dissipation},
\newblock \bibinfo{journal}{Journal of Scientific Computing}
  \bibinfo{volume}{29} (\bibinfo{year}{2006}) \bibinfo{pages}{375--404}.
\bibitem[{Zienkiewicz and Taylor(2000)}]{Zien2000}
\bibinfo{author}{O.~Zienkiewicz}, \bibinfo{author}{R.~Taylor},
\newblock \bibinfo{title}{Finite element method: Volume 1},
\newblock in: \bibinfo{booktitle}{Butterworth-Heinemann; 5th edition}.
\bibitem[{RA and CR(1990)}]{Horn1990}
\bibinfo{author}{H.~RA}, \bibinfo{author}{J.~CR},
\newblock \bibinfo{title}{Matrix analysis, cambridge university press},
\newblock in: \bibinfo{booktitle}{Cambridge, UK}.
\bibitem[{Frenander and Nordstr\"om(2018)}]{Frenander2018}
\bibinfo{author}{H.~Frenander}, \bibinfo{author}{J.~Nordstr\"om},
\newblock \bibinfo{title}{Spurious solutions for the advection-diffusion
  equation using wide stencils for approximating the second derivative},
\newblock \bibinfo{journal}{Numerical Methods for PDEs} \bibinfo{volume}{34}
  (\bibinfo{year}{2018}) \bibinfo{pages}{501--517}.
\bibitem[{Malan et~al.(2002)Malan, Lewis, and Nithiarasu}]{Malan2002}
\bibinfo{author}{A.~Malan}, \bibinfo{author}{R.~Lewis},
  \bibinfo{author}{P.~Nithiarasu},
\newblock \bibinfo{title}{An improved unsteady, unstructured, artificial
  compressibility, finite volume scheme for viscous incompressible flows: Part
  i. theory and implementation},
\newblock \bibinfo{journal}{International Journal for Numerical Methods in
  Engineering} \bibinfo{volume}{54} (\bibinfo{year}{2002})
  \bibinfo{pages}{695--714}.
\bibitem[{Eriksson et~al.(2011)Eriksson, Abbas, and Nordstr\"om}]{Sofia2011}
\bibinfo{author}{S.~Eriksson}, \bibinfo{author}{Q.~Abbas},
  \bibinfo{author}{J.~Nordstr\"om},
\newblock \bibinfo{title}{A stable and conservative method for locally adapting
  the design order of finite difference schemes},
\newblock \bibinfo{journal}{Journal of Computational Physics}
  \bibinfo{volume}{230} (\bibinfo{year}{2011}) \bibinfo{pages}{4216--4231}.
\bibitem[{Jiang and Shu(1996)}]{Weno3_1996}
\bibinfo{author}{G.~Jiang}, \bibinfo{author}{C.~Shu},
\newblock \bibinfo{title}{Efficient implementation of weighted eno schemes},
\newblock \bibinfo{journal}{Journal of Computational Physics}
  \bibinfo{volume}{126} (\bibinfo{year}{1996}) \bibinfo{pages}{202--228}.
\bibitem[{Fisher et~al.(2013)Fisher, Carpenter, Nordstr\"om, Yamaleev, and
  Swanson}]{JanNonLin2013}
\bibinfo{author}{T.~C. Fisher}, \bibinfo{author}{M.~H. Carpenter},
  \bibinfo{author}{J.~Nordstr\"om}, \bibinfo{author}{N.~K. Yamaleev},
  \bibinfo{author}{C.~Swanson},
\newblock \bibinfo{title}{Discretely conservative finite-difference
  formulations for nonlinear conservation laws in split form: Theory and
  boundary conditions},
\newblock \bibinfo{journal}{Journal of Computational Physics}
  \bibinfo{volume}{234} (\bibinfo{year}{2013}) \bibinfo{pages}{353--375}.

\end{thebibliography}

\pagebreak 

\end{document}